\documentclass[preprint,12pt]{elsarticle}
\biboptions{sort&compress}

\usepackage{amssymb}
\usepackage{amsmath}
\usepackage{xcolor}
\usepackage{lmodern}
\usepackage{import}
\usepackage{booktabs,subcaption}

\journal{Finite Elements in Analysis and Design}
\newcommand{\mbn}[1]{\boldsymbol{\mathrm{#1}}}

\begin{document}

\begin{frontmatter}


 \ead{wolfgang.dornisch@mv.rptu.de}
 \cortext[cor1]{Current affiliation: Rheinland-Pfälzische Technische Universität Kaiserslautern-Landau, Gottlieb-Daimler-Str., 67663 Kaiserslautern, Germany}
  \cortext[cor2]{Corresponding author}

\title{A rotation-based geometrically nonlinear spectral Reissner--Mindlin shell element}

 \author[label1,label2]{Nima Azizi\corref{cor1}}
 \affiliation[label1]{organization={Rheinland-Pfälzische Technische Universität Kaiserslautern-Landau},
             addressline={Gottlieb-Daimler-Str.},
             city={Kaiserslautern},
             postcode={67663},
             country={Germany}}

 \affiliation[label2]{organization={Brandenburgische Technische Universität Cottbus-Senftenberg},
             addressline={Konrad-Wachsmann-Allee~2},
             city={Cottbus},
             postcode={03046},
             country={Germany}}

 \author[label1]{Wolfgang Dornisch\corref{cor2}}


\begin{abstract}
In this paper, we propose a geometrically nonlinear spectral shell element based on Reissner--Mindlin kinematics using a rotation-based formulation with additive update of the discrete nodal rotation vector. The formulation is provided in matrix notation in detail. The use of a director vector, as opposed to multi-parameter shell models, significantly reduces the computational cost by minimizing the number of degrees of freedom. Additionally, we highlight the advantages of the spectral element method (SEM) in combination with Gauss-Lobatto-Legendre quadrature regarding the computational costs to generate the element stiffness matrix. To assess the performance of the new formulation for large deformation analysis, we compare it to three other numerical methods. One of these methods is a non-isoparametric SEM shell using the geometry definition of isogeometric analysis (IGA), while the other two are IGA shell formulations which differ in the rotation interpolation. All formulations base on Rodrigues' rotation tensor.  Through the solution of various challenging numerical examples, it is demonstrated that although IGA benefits from an exact geometric representation, its influence on solution accuracy is less significant than that of shape function characteristics and rotational formulations. Furthermore, we show that the proposed SEM shell, despite its simpler rotational formulation, can produce results comparable to the most accurate and complex version of IGA. Finally, we discuss the optimal SEM strategy, emphasizing the effectiveness of employing coarser meshes with higher-order elements.
\end{abstract}




\newtheorem{remark}{Remark}

\end{frontmatter}



\section{Introduction}
Due to the widespread usage and importance of shell structures, both in nature and industry, the study of shell structures has remained an active area of research. 
Although the finite element method (FEM) has been and continues to be widely used in the numerical simulation of shell structures, more recent and versatile methods, such as isogeometric analysis (IGA)~\cite{hughes_2005}, are also being developed and utilized in shell analysis as well as in many other applications~\cite{nguyen2015}. In most cases, Non-Uniform Rational B-splines (NURBS) are used within IGA~\cite{Piegl.1997}. IGA allows for a straightforward implementation of Kirchhoff-Love shell elements~\cite{Kiendl.2009}. A general drawback of  Kirchhoff-Love elements is their inability to capture accurate results in regions of high plastic deformation~\cite{Benson.2013}. Thus, an isogeometric Reissner-Mindlin element was developed in~\cite{Benson.2010}, where the use of high-order elements is proposed as a remedy for locking effects. However,~\cite{Oesterle.2017, zou2020, zou2022} emphasize that despite using high-order elements, locking effects persist, and various methods are presented to eliminate them. Other methods are also proposed to alleviate membrane locking in isogeometric Kirchhoff-Love elements, which inherently are free of shear locking and suffer only from membrane locking~\cite{ZOU2021, casquero2023, leonetti2024}. Focusing on director-based isogeometric Reissner-Mindlin shell elements, two controversial points remain. The first one is the calculation of the nodal coordinate systems at the control point based on coordinate systems at integration points, which can be determined to be exactly tangential~\cite{Dornisch.2013}. The second, and more important point is the method of calculation of the rotated director vector at the integration points. While~\cite{Simo.1989b} demonstrates that the simple interpolation of the rotated director vector using Rodrigues' formula is sufficient for lower-order Lagrange-based finite elements, the large support of IGA NURBS shape functions violates this condition and may introduce errors in the analysis in the form of artificial thinning \cite{Dornisch.2013,Dornisch.2015b}. To eliminate this effect, one solution is to calculate the rotated director vectors based on continuous rotation interpolation \cite{Dornisch.2013}.
\par If a method can address the difficulty of accurately calculating director vectors at integration points without adding more complexity and, as a potential initial remedy for locking, is capable of implementing very high-order elements sufficient to eliminate the locking effect \cite{payette_2014}, it reduces computational costs and potentially simplifies the formulation. A natural choice is the spectral element method, in which the element nodes coincide with integration points~\cite{patera_1984, poz_2005}, and the location of these points are calculated
based on the root of different polynomials. There are various types of spectral element methods (SEM), but we focus on the use of Legendre nodes alongside Gauss-Lobatto-Legendre (GLL) quadrature. This method allows for the use of very high-order shape functions while still achieving stable and accurate results~\cite{PITTON2018440}. In~\cite{nima_2024} we proposed a linear Reissner--Mindlin spectral shell element based on a NURBS geometry description. Notably, the use of very high orders in SEM significantly resolves the shear locking problem \cite{nima_2024} and yields stable computations up to an order of around $p=20$. On the one hand, SEM was used before to analyze shell structures in~\cite{payette_2014, Gutierrez2016, zak2018}. However, in these studies, attempts to incorporate thickness stretch, transverse shear deformation and higher-order shear deformation -- in order to prevent locking and achieve a more comprehensive representation of deformations -- result in the automatic exclusion of the concept of director vector. These points are achieved through multi-parameter expansion of the deformation, which, consequently, significantly increases the number of degrees of freedom in the system of equations. On the other hand, several publications about Reissner--Mindlin plates have been presented in the literature, e.g.~\cite{Cohen.2007,Ambati.2025}. However, a spectral version of the geometrically nonlinear, rotation-based Reissner--Mindlin shell formulation, as proposed by~\cite{Simo.1989,Simo.1989b} for the low-order Lagrange-based finite element method, is still outstanding.
\par In this contribution, we are proposing a geometrically nonlinear shell formulation based on Reissner--Mindlin kinematics and an additive update of the director vector at the nodes using the full Rodrigues' rotation tensor~\cite{argyris_1982}. This represents the first geometrically nonlinear Reissner--Mindlin SEM shell formulation. In contrast to earlier works on SEM shells~\cite{payette_2014, Gutierrez2016, zak2018}, we use only five degrees of freedom (DOF) per node to obtain fast computations. The primary focus of using high orders of basis functions is not high preciseness as in the higher order shear-deformable SEM shell formulations mentioned above, but the alleviation of locking effects.  In contrast to our earlier work~\cite{nima_2024}, we do not use NURBS basis functions within the formulation, but stick to a purely isoparametric formulation. However, we compute the location of the nodes from the exact NURBS geometry. Since in SEM integration points and nodes coincide, the physical location of the integration points in the reference configuration is always exact. But in the formulation, the derivatives of physical location and director vector occur. Since we interpolate these with Lagrange basis functions, the derivatives of these quantities are approximated. This contribution will detail the fundamental assumptions and formulations employed in the large deformation analysis of the proposed isoparametric spectral Reissner-Mindlin shell formulation in the first part of Sec.~\ref{sec_2}. Subsequently, to provide a proper basis for the numerical examples, we will shortly present and explain a different variant of SEM based on our previous work~\cite{nima_2024}, and two IGA methods based on~\cite{Dornisch.2016}, which differ in the rotation interpolation. In Sec.~\ref{sec_3}, we will propose a strategy to speed up SEM computations and compare the efficiency of the methods from a theoretical perspective by counting required multiplications. Various challenging numerical example are solved in Sec.~\ref{sec_4}, where the examples are chosen to have rising complexity by means of high and highly changing curvature. Primarily, we want show the accuracy and efficiency of the newly proposed shell formulation. As a side effect of the comparison with the IGA shell formulations, we also draw conclusions regarding the importance of the exact geometry representation in IGA. We close with Discussion and Conclusion in Sec.~\ref{sec_conclusion}.

\section{Theoretical background and formulation}
\label{sec_2}
Before explaining in detail the discretized form of the proposed isoparametric spectral Reissner--Mindlin shell formulation, it is helpful to provide a general overview of the underlying shell theory. Therefore, we present the continuous equations in Sec.~\ref{subsec.2-1}, followed by a detailed elaboration of the discretized equations and the resulting matrices in Sec. \ref{subsec.2_2}.

\subsection{Reissner--Mindlin shell element formulation}
\label{subsec.2-1}

\subsubsection{Kinematics}
\label{subsec.2-1-1}
\begin{figure}[ht] 
	\center{\includegraphics[width=\textwidth]{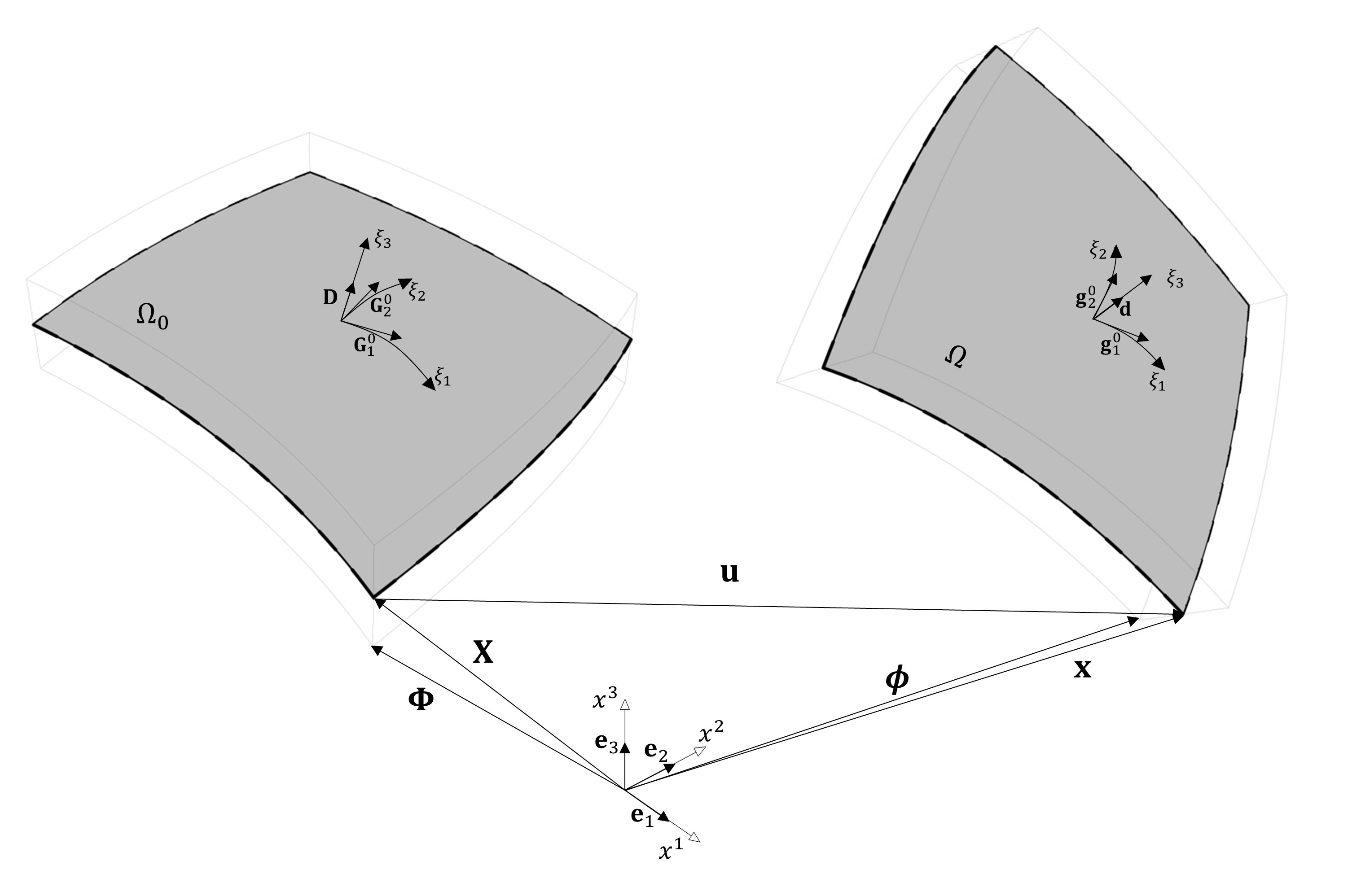} }
	\caption{\centering Undeformed (reference) and deformed (current) configuration of the shell}\label{fig_shell_conf}
\end{figure}

As depicted in Fig.~\ref{fig_shell_conf}, the mathematical description of an undeformed configuration of a shell body, $\mathrm{B}_0$, can be formulated by its undeformed midsurface, $\mathrm{\Omega_0}
$, and a unit director vector, $\mbn{D}$, which is usually perpendicular to the midsurface in the undeformed configuration. So, the location vector in undeformed configuration, $\mbn{\Phi}$, of any material point inside the shell body  is defined as below

\begin{equation}\label{shell_conf}
	\mbn{\Phi}(\xi^i) = \Phi^i \mbn{e}_i= \mbn{X}(\xi^\alpha)+\xi^3 \mbn{D}(\xi^\alpha)\,,
\end{equation}
where  $\Phi^i$ and $\mbn{e}_i$ are Cartesian components of $\mbn{\Phi}$ and the Cartesian unit base vectors, respectively. $\mbn{X}$ denotes the position vector of a point on the midsurface of the shell. $\xi^{i}$ is the coordinate of a material point in the curvilinear covected coordinate system. The convected coordinate system can be considered as the natural or parametric coordinate system of a finite element formulation. If not specified differently, Greek indices range from 1 to 2 and Latin indices from 1 to 3. $\xi^3$ is the parametric coordinate in the direction of thickness of the shell surface and runs from $-h/2$ to $h/2 $, where $h$ is the thickness of the shell. For simplicity, we focus on shell structures with constant thickness, but variations in thickness can be easily incorporated into the formulation. The director vector is calculated by
\begin{equation}\label{calculation_dir}
	\mbn{D}(\xi^\alpha) = \frac{\mbn{G}^0_1(\xi^\alpha) \times \mbn{G}^0_2(\xi^\alpha)}{ 	\left| \mbn{G}^0_1(\xi^\alpha) \times \mbn{G}^0_2(\xi^\alpha)	\right|}\,.
\end{equation}
Accordingly, the director vector has unit length, i.e.,
\begin{equation}\label{director}
	\left| \mbn{D}(\xi^\alpha) \right| = 1\,.
\end{equation}
$\mbn{G}^0_\alpha$ is the covariant tangent vector to the undeformed  midsurface of the shell and defined as
\begin{equation}\label{G_0}
	\mbn{G}^0_\alpha = \frac{\partial\mbn{X}}{\partial\xi^\alpha}\,,
\end{equation}
and the covariant vector at any point through the thickness of the shell is
\begin{equation}\label{G}
	\mbn{G}_\alpha = \frac{\partial\mbn{\Phi}}{\partial\xi^\alpha}=	\mbn{G}^0_\alpha+\xi^3 \mbn{D},_{\xi^\alpha} \quad \mathrm{and} \quad \mbn{G}_3=\mbn{D}\,.
\end{equation}
\par As depicted in Fig.~\ref{fig_shell_conf}, the position vector of any point within the shell body in a deformed configuration with the deformed mid-surface, $\mathrm{\Omega}$, is 
 \begin{equation}\label{shell_conf_deform}
     \mbn{\phi}(\xi^i) = \mbn{x}(\xi^\alpha)+\xi^3 \mbn{d}(\xi^\alpha)\,.
 \end{equation}
Here $\mbn{x}$ denotes the position vector of a point on the midsurface in deformed configuration and $\mbn{d}$ is the rotated director vector. Consequently, the displacement vector of the midsurface  $\mbn{u}$ is given by
\begin{equation}\label{u}
	\mbn{u}(\xi^{\alpha})=\mbn{x}(\xi^{\alpha})-\mbn{X}(\xi^{\alpha})\,.
\end{equation}
In Reissner--Mindlin shell theory, the first order shear deformation theory is utilized, so the director is not necessarily perpendicular to the midsurface in the deformed configuration. Also, in this theory, there is no extension in the thickness of the shell and the domain of $\xi^3$ always ranges from $-h/2$ to $h/2 $, therefore, the following equation must always hold
\begin{equation}\label{director_t}
	\left| \mbn{d}(\xi^\alpha) \right| = 1\,.
\end{equation}

\noindent Similar to Eq.~\eqref{G}, the covariant tangent vector in deformed configuration to midsurface shown in Fig.~\ref{fig_shell_conf} is
\begin{equation}\label{g_0}
	\mbn{g}^0_\alpha = \frac{\partial\mbn{x}}{\partial\xi^\alpha}
\end{equation}
and the covariant vector at any point through the thickness of the shell is
\begin{equation}\label{g}
	\mbn{g}_\alpha = \frac{\partial\mbn{\phi}}{\partial\xi^\alpha}=	\mbn{g}^0_\alpha+\xi^3 \mbn{d},_{\xi^\alpha}\quad \mathrm{and} \quad \mbn{g}_3=\mbn{d}\,.
\end{equation}
Obviously, in the presence of shear deformation, it is not possible to use covariant tangent vectors to generate director vector in deformed configuration and the method of calculation of director vectors in this configuration will be discussed in Sec. \ref{subsec.2_2}.

\subsubsection{General strain formulation}
\label{subsec.2-1-2}
The Green-Lagrange strain tensor is chosen to analyze the large deformation and large rotation of shells. If we define $\mbn{F}=\mathrm{Grad}\;\mbn{\phi}$, then the Green-Lagrange strain tensor is
\begin{equation}\label{green_strain}
	\mbn{E} = \frac{1}{2}(\mbn{F}^T \mbn{F}-\mbn{1}),
\end{equation}
where $\mbn{1}$ is the unit tensor, which takes the familiar form of a 2-dimensional matrix with elements on the main diagonal equal to 1 in the Cartesian coordinate system. Regarding the point that the unit tensor in Eq.~\eqref{green_strain} is indeed the metric tensor, and using the contravariant basis vectors $\mbn{G}^i$, the strain components in the curvilinear coordinate system are \cite{wagner_2005}
\begin{equation}\label{green_cov}
	\mbn{E}=E_{ij} \mbn{G}^i \otimes \mbn{G}^j,\quad E_{ij}=(\mbn{\phi}_{,\xi^i} \cdot \mbn{\phi}_{,\xi^j}-\mbn{\Phi}_{,\xi^i} \cdot \mbn{\Phi}_{,\xi^j})\,.
\end{equation}
Here, comma shows the partial derivative with respect to curvilinear coordinate $\xi^i$ and the contravariant base vector is
\begin{equation}\label{G-contra}
	\mbn{G}^\alpha = \frac{\partial\xi^\alpha}{\partial\mbn{\phi}}\,.
\end{equation}
Using Eqs.~\eqref{shell_conf} and \eqref{shell_conf_deform} in Eq.~\eqref{green_cov} and neglecting higher order curvature yields
\begin{equation} 
	\begin{split}\label{shell_strain_1}
	& E_{\alpha \beta}= \epsilon_{\alpha \beta}+\xi^3 \kappa_{\alpha \beta}  \\
    & 2E_{3 \alpha}=\gamma_{\alpha} \\
	& E_{33}=0 \,,
	\end{split}
\end{equation}
where the shell strains are
\begin{equation} \label{shell_strain_2}
	\begin{split}
		 &\epsilon_{\alpha \beta}= \frac{1}{2}(\mbn{x}_{,\xi^\alpha}\cdot \mbn{x}_{,\xi^\beta}-\mbn{X}_{,\xi^\alpha} \cdot \mbn{X}_{,\xi^\beta})  \\
		& \kappa_{\alpha \beta}=\frac{1}{2}(\mbn{x}_{,\xi^\alpha} \cdot \mbn{d}_{,\xi^\beta}+\mbn{x}_{,\xi^\beta} \cdot \mbn{d}_{,\xi^\alpha}-(\mbn{X}_{,\xi^\alpha} \cdot \mbn{D}_{,\xi^\beta}+\mbn{X}_{,\xi^\beta} \cdot \mbn{D}_{,\xi^\alpha})) \\
		&\gamma_{\alpha}=\mbn{x}_{,\xi^\alpha} \cdot \mbn{d}-\mbn{X}_{,\xi^\alpha} \cdot \mbn{D} \\
	\end{split}
\end{equation}
with $\epsilon_{\alpha \beta}$, $\kappa_{\alpha \beta}$ and $\gamma_{\alpha}$ being the membrane, curvature and shear strains, respectively.

\subsection{Various discretization formulations}
\label{subsec.2_2}
The main contribution of this paper is to propose an isoparametric spectral Reissner-Mindlin shell formulation based on a rotational formulation for the director vector. We present this formulation in full detail, followed by a short summary of the main features of the other formulations used in the numerical examples. 
\subsubsection{Spectral Element Method, Isoparametric formulation (SEMI)}\label{semi}
\begin{figure}[ht] 
	\center{\includegraphics[scale=0.7]{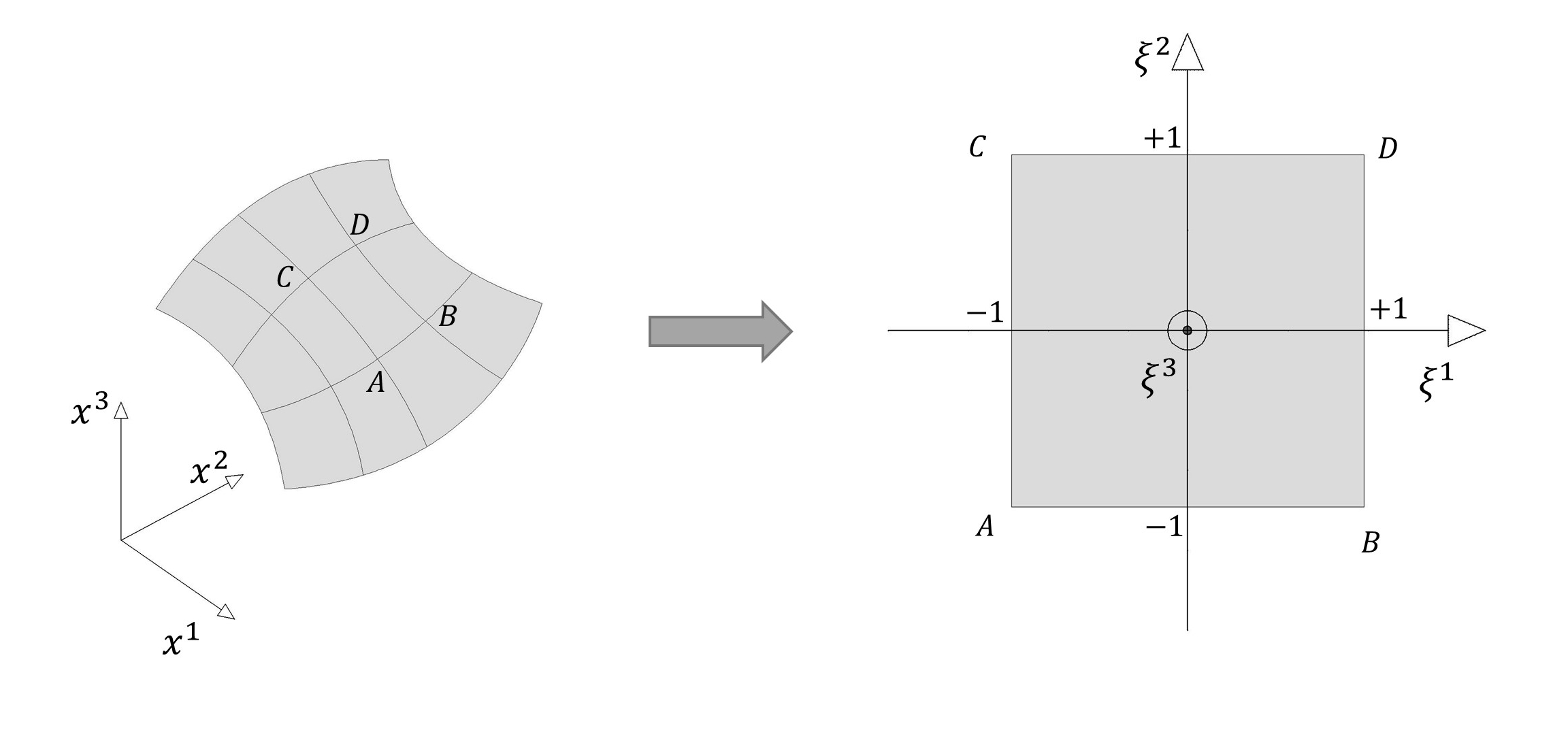} }
	\caption{\centering Physical space and SEMI parametric space }\label{fig_map_semi}
\end{figure}

This contribution bases on the isoparametric Legendre spectral finite element method. Lagrange shape functions are used to interpolate unknown fields and the coordinates of  unevenly spaced nodes of elements are calculated based on  the roots of Legendre polynomials \cite{patera_1984,poz_2005}. If we take the number of interpolating nodes as $(p+1)$, where \textit{p} is the order of element, then the $p$th order 1D Lagrange shape function for the \textit{i}th interpolating point is calculated in natural (parametric) coordinate system according to 
\begin{equation}\label{lagrange}
	l^{p}_i(\xi) = \frac{\prod_{j=1, j\ne i}^{p+1}(\xi-\xi^j)}{\prod_{j=1, j\ne i}^{p+1}(\xi^i-\xi^j)}     \quad     -1\leq \xi\leq+1 \quad \mathrm{and} \qquad 1\leq i\leq p+1\,,
\end{equation}
where $\xi^j$ is the coordinate of the \textit{j}th node in the SEM parametric space. After some algebra, the differentiation of Eq.~(\ref{lagrange}) yields the derivatives
\begin{equation}\label{lagrange_der}
	\frac{\mathrm{d} l^{p}_i(\xi)}{\mathrm{d}\xi}=\frac{\sum_{j=1, j\ne i}^{p+1}\prod_{k=1, k\ne i, k\ne j}^{p+1}(\xi-\xi^k)}{\prod_{j=1, j\ne i}^{p+1}(\xi^i-\xi^j)} \quad 1\leq i\leq p+1
\end{equation}
of the Lagrange shape functions. Details about the calculations can be found in \cite{poz_2005}.
Within our SEM formulation, we use Gauss-Lobatto-Legendre (GLL) numerical integration, for which
\begin{equation}
	\int_{-1}^{+1} f(\xi) \mathrm{d}\xi \simeq \sum_{i=1}^{p+1}f(\xi^i)w_i   \qquad 1\leq i\leq p+1  
\end{equation}
with 
\begin{equation}
	w_i=\int_{-1}^{+1}l^{p}_i(\xi)\mathrm{d}\xi
\end{equation}
holds. The 2D Lagrange shape functions can be calculated by the tensor product of vectors of 1D shape functions from Eq.~\eqref{lagrange} as follows
\begin{equation}\label{lagr_2d}
    N_I(\xi^1,\xi^2)= l_i(\xi^1) l_j(\xi^2) \qquad 1\leq I \leq n_{en} \,,
\end{equation} 
where $ n_{en}=(p+1)^2$ represents the total number of nodes in an element, assuming that the element has the same order $p$ in both parametric directions. As shown in Fig.~\ref{fig_map_semi}, the geometry of each element $ABCD$ is mapped from physical space to parametric space by utilizing these 2D shape functions.
\par
\begin{figure}[ht] 
	\center{\includegraphics[scale=0.45]{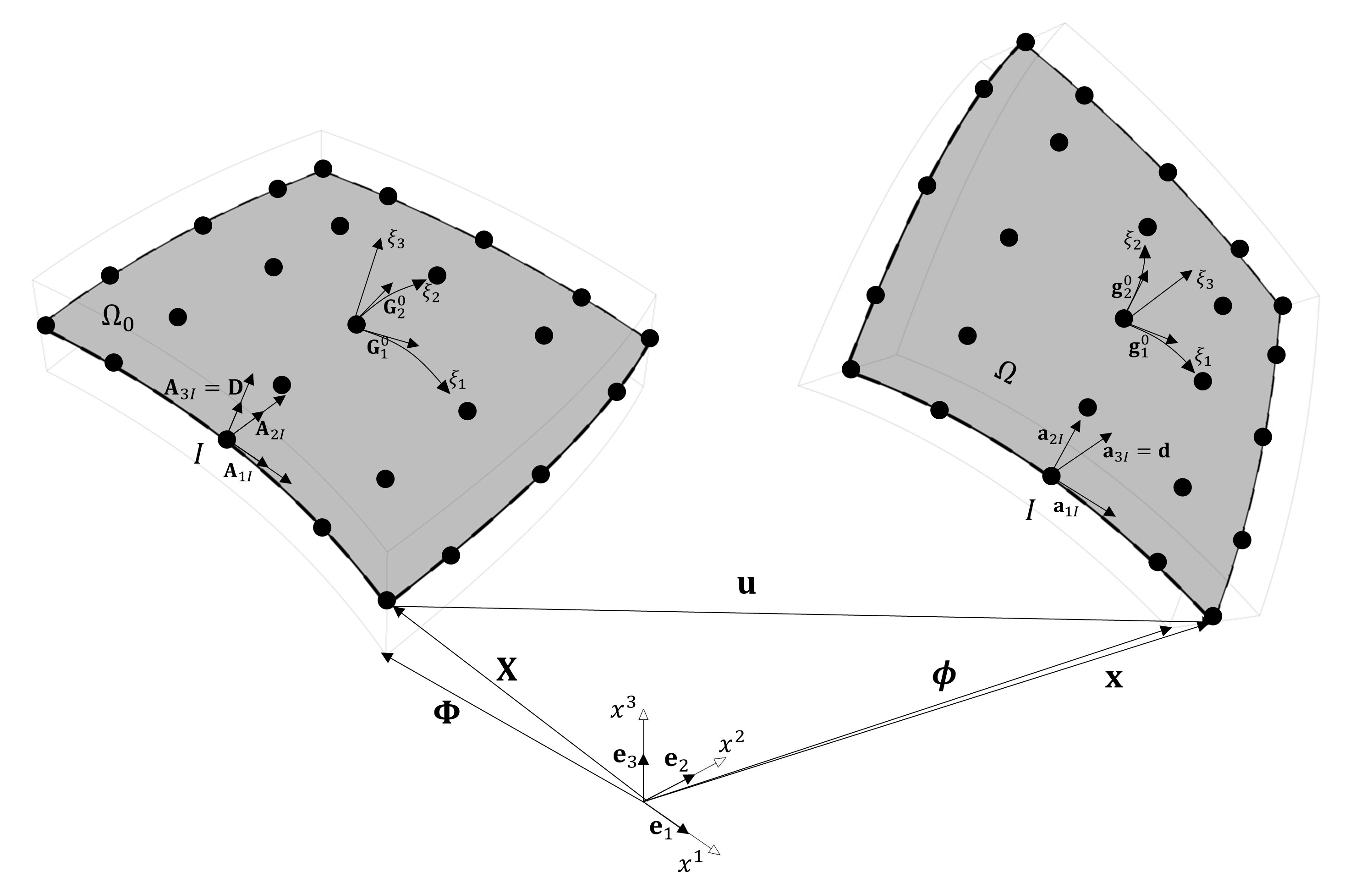} }
	\caption{\centering Undeformed (reference) and deformed (current) configuration of the shell}\label{fig_SEMI_element}
\end{figure}

 To avoid any ambiguity, it is essential to first clarify the coordinate systems and degrees of freedom (DOF) in greater detail before deriving the stiffness matrices and explaining the solution process. With these foundational concepts established, the subsequent equations will be easier to follow. Because the Reissner--Mindlin shell theory is not capable of generating drilling stiffness, in all of the discretization methods investigated in this paper, at each node on the surface of the shell the drilling DOF is neglected. To easily exclude the drilling DOF, a local Cartesian coordinate system with the third axis perpendicular to the midsurface is defined at each node  as depicted for node $I$ in Fig.~\ref{fig_SEMI_element}.  In this figure, for demonstration purposes, only one 4th-order element has been used for modeling the shell. By neglecting the component of a rotation vector around the axis 3 in the local nodal coordinate system, drilling will be excluded. These local nodal Cartesian systems also should be constructed for each node in a way that the rotational boundary condition can be applied easily. We follow the method described in \cite{Bathe.2006} and presented in \ref{app1} to create these systems for each node and calculate the unit base system $\mbn{A}_{I}=[\mbn{A}_{1I}, \mbn{A}_{2I}, \mbn{A}_{3I}]$ with the three unit base vectors $\mbn{A}_{iI}$. This local nodal system rotates during the deformation, and in deformed configuration, it is represented by  $[\mbn{a}_{1I}, \mbn{a}_{2I}, \mbn{a}_{3I}]$ being calculated by using Rodrigues' orthogonal transformation~\cite{argyris_1982}. The Rodrigues' transformation, $\mbn{R}_I$, in the form which is used in our contribution is
\begin{equation} \label{rodrig}
  \begin{split}
	\mbn{R}_I&=\mbn{1}+\frac{\mathrm{sin}\, \omega_I}{\omega_I} \mbn{\Omega}_I+\frac{1-\mathrm{cos}\, \omega_I}{\omega_I^2} \mbn{\Omega}_I^2 \\
	\mbn{\Omega}_I&=
	\begin{bmatrix}
		&0 &-\omega_{zI} &\omega_{yI} \\
		&\omega_{zI} &0 &-\omega_{xI} \\
		&-\omega_{yI} &\omega_{xI} &0 
	\end{bmatrix} \\
	\omega_I&= | \mbn{\omega}_I |\,.
  \end{split}
\end{equation}
Here $\mbn{\omega}_I$ represents the total rotation vector experienced by each point $I$ during the deforming from $\Omega_0$  to $ \Omega$ configuration and $\omega_{xI}$, $\omega_{yI}$ and $\omega_{zI}$ are the global Cartesian components of $\mbn{\omega}_I$. We use $x$, $y$ and $z$ symbols to represent the components of a vector in global Cartesian coordinate and to distinguish them from the local nodal Cartesian components. There are some subtle points about the singularity in Eq.~\eqref{rodrig}, which are explained in detail in \cite{wagner_2005} and we refrain from repeating them. Then, the rotated local nodal Cartesian coordinate in deformed configuration is \cite{wagner_2005}
\begin{equation}\label{a_A}
	\mbn{a}_{iI} = \mbn{R}_I  \mbn{A}_{iI}\,.
\end{equation}

\par
\begin{figure}[ht] 
	\center{\includegraphics[scale=0.7]{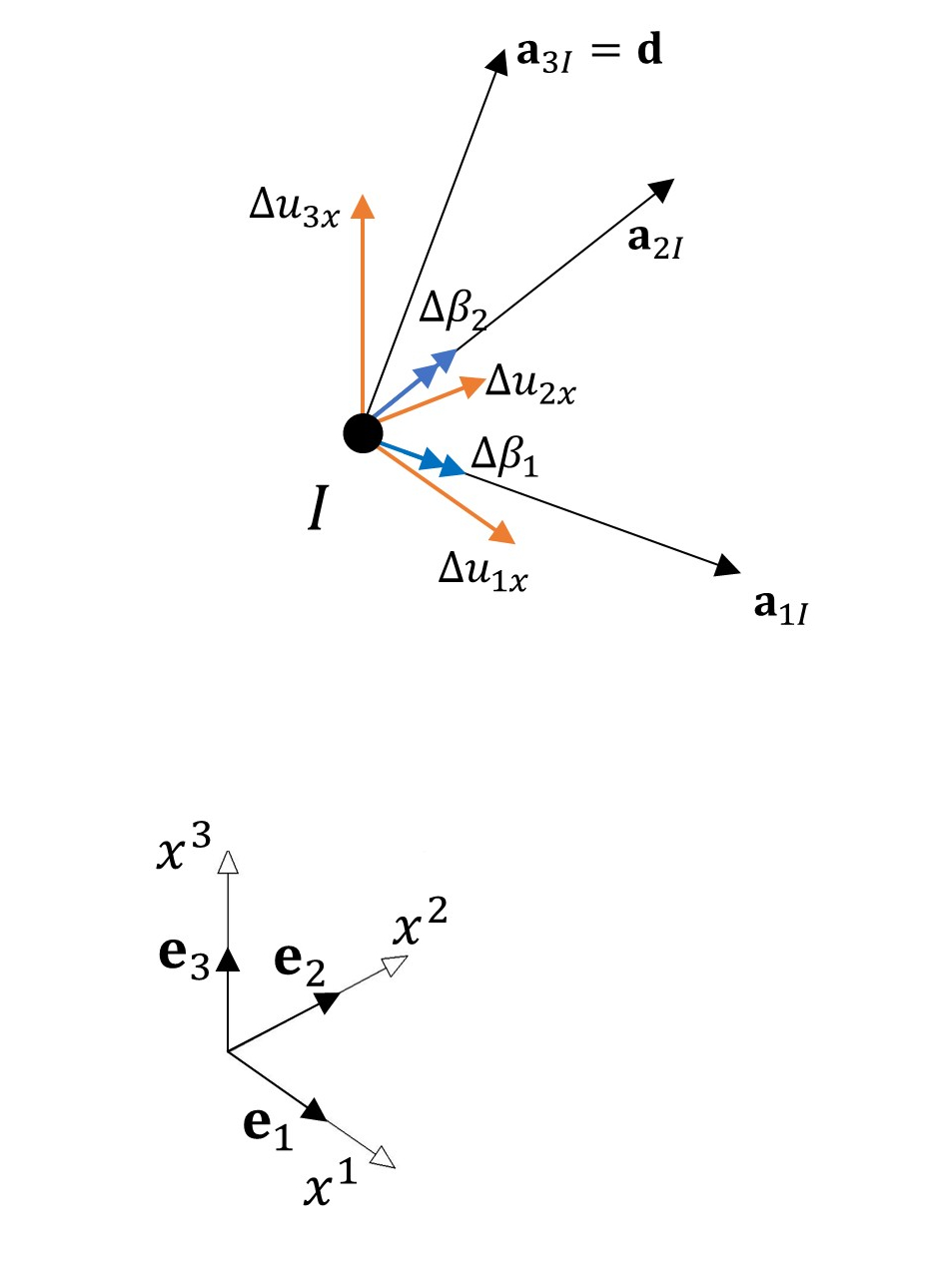} }
	\caption{\centering DOFs of node $I$}\label{fig_DOFs}
\end{figure}

 Neglecting drilling, each node has 5 DOFs unless the node lies on the intersection of two shells \cite{wagner_2005}, which will be clarified later. These DOFs consist of three translational and two rotational DOFs. These two rotational DOFs, which are defined in the local nodal coordinate system, form a vector at each node $I$, which is denoted as $\Delta\mbn{\beta}_I$. The translational DOFs, on the other hand, are the components of the  displacement vector in global Cartesian coordinates. So, for a node outside of the shell intersection the increment of DOFs is
\begin{equation}
	\Delta\mbn{v}_I=\begin{bmatrix}
		\Delta\mbn{u}\\
		\Delta\mbn{\beta}
	\end{bmatrix}_I=
	\begin{bmatrix}
		\Delta u_{xI} \\
		\Delta u_{yI} \\
		\Delta u_{zI} \\
		\Delta\beta_{1I} \\
		\Delta\beta_{2I} \\
	\end{bmatrix}\,,
\end{equation}

\noindent which are depicted in Fig.~\ref{fig_DOFs}. Since the vector $\mbn{a}_{iI}$ varies at each instant of time, we cannot directly use $\Delta\mbn{\beta}_I$ to additively update the total rotation vector, $\mbn{\omega}_I$, at each node. Therefore, it is necessary to transform $\Delta\mbn{\beta}_I$ into the constant global Cartesian system as follows \cite{wagner_2005}

\begin{equation}\label{beta_to_omega}
	\Delta\mbn{\omega}_I= \mbn{T}_{3I}\Delta\mbn{\beta}_I \,,
\end{equation}
where
\begin{equation}\label{eq_T3I}
	\mbn{T}_{3I}=
	\begin{cases}
		\mbn{1}_{3\times 3} & \text{for nodes on shell intersections} \\
		[\mbn{a}_{1I}, \mbn{a}_{2I}]_{(3 \times 2)} & \text{for all other nodes}
	\end{cases}
\end{equation}
and
\begin{equation}
	\Delta\mbn{\beta}_I=
	\begin{cases}
		[\Delta\beta_{xI}, \Delta\beta_{yI}, \Delta\beta_{zI}]^\mathrm{T} &\text{for nodes on shell intersections}\\
		[\Delta\beta_{1I}, \Delta\beta_{2I}]^\mathrm{T} & \text{for all other nodes}
	\end{cases}\,.
\end{equation}
Because at shell intersections each side can provide stiffness for the drilling of the other side, drilling must not be excluded from the calculations. Therefore, the three rotational DOFs are considered for nodes at the intersections. Finally the displacement and rotation vectors for node $I$ at the end of iteration $k$ are 

\begin{subequations}
	\begin{align}
	\mbn{u}^{(k)}_I=\mbn{u}_I^{(k-1)}+\Delta\mbn{u}^{(k)}_I\\
	\mbn{\omega}^{(k)}_I = \mbn{\omega}^{(k-1)}_I+\Delta\mbn{\omega}^{(k)}_I
	\end{align}
\end{subequations}
\begin{remark} 
It is noteworthy to mention that, generally, the components of a rotation vector do not represent rotations around the axes \cite{argyris_1982} unless the rotation is small enough. This is the assumption that we make for each Newton--Raphson iteration step during the solution process. So, zero rotation around aixs $ \mbn{a}_3$ is approximately equivalent to setting the third component of rotation in the nodal local Cartesian system to 0 and hence excluding the drilling.
\end{remark}
\par Using the isoparametric assumptions and regarding the explained DOFs for each node, the required interpolated functions in each element are
\begin{subequations}\label{interpolation}
	\begin{align}
		&\mbn{X}^h = \sum_{I=1}^{n_{en}} N_I\,\mbn{X}_I\, ,\\ 
		&\mbn{x}^h =\sum_{I=1}^{n_{en}} N_I \, \mbn{x}_I \, ,\\
		&\mbn{d}^h =\sum_{I=1}^{n_{en}} N_I \, \mbn{d}_I \, , \label{d_inter}\\[0.1cm]
		&\mbn{d}_I =\mbn{R}_I \, \mbn{D}_I\,. \label{d_rot}
	\end{align}
\end{subequations}
The derivatives of the 2D lagrange shape functions are also obtained based on Eq.~\eqref{lagrange_der} with respect to the convective bases $\xi^\alpha$. Eq.~\eqref{d_rot} is based on Eq.~\eqref{a_A} and $\mbn{x}_I$ are calculated based on Eq.~\eqref{u}. The derivatives of the interpolated functions with respect to convected coordinates, $\xi^\alpha$, are
\begin{subequations}\label{interpolation_der}
	\begin{align}
		&\mbn{X}^h_{,\xi^\alpha} = \sum_{I=1}^{n_{en}} N_{I,\xi^{\alpha}}\,\mbn{X}_I \, ,\\ 
		&\mbn{x}^h_{,\xi^\alpha} =\sum_{I=1}^{n_{en}} N_{I,\xi^\alpha} \, \mbn{x}_I \, ,\\
		&\mbn{d}^h_{,\xi^\alpha} =\sum_{I=1}^{n_{en}} N_{I,\xi^\alpha} \, \mbn{d}_I \, \label{d_inter_der}.
	\end{align}
\end{subequations}

\par To establish the constitutive equations in a simple manner, it is common practice to use local Cartesian basis systems for the calculation of strains. Depending on the intended material law, different choices are possible. The only requirement is that the plane spanned by $\mbn{A}_{1}$ and $\mbn{A}_{2}$ must be tangential to the shell surface in the considered location. Since in SEMI nodes and integration points coincide, we can directly use the local nodal coordinate systems $\mbn{A}_I$, of which the construction is defined in \ref{app1}. It is to be noted that in cases of anisotropic material, it might be advantageous to choose local basis systems which are aligned to the material directions and can differ from $\mbn{A}_I$, which is possible for the proposed formulation in a straightforward manner by replacing $\mbn{A}_{iI}$ in Eq.~\eqref{eq:jacob} with the basis vectors of another system. So, instead of Eqs.~\eqref{shell_strain_1} and \eqref{shell_strain_2}, which are presented in a general curvilinear coordinate system, the Green-Lagrange strains will now be calculated in the local nodal Cartesian coordinate systems~$\mbn{A}_I$. Consequently, derivatives of the position vectors will be needed with respect to the coordinates of these systems. Therefore, the Jacobian matrix is needed to map from the general natural coordinate systems to the local nodal Cartesian systems. It is given by \cite{wagner_2005}
\begin{equation}
	\mbn{J}=
	\begin{bmatrix} 
		\mbn{X}_{,\xi^1} \cdot \mbn{A}_{1I} & \mbn{X}_{,\xi^1} \cdot \mbn{A}_{2I}\\
		\mbn{X}_{,\xi^2} \cdot \mbn{A}_{1I} & \mbn{X}_{,\xi^2} \cdot \mbn{A}_{2I}
	\end{bmatrix}.
    \label{eq:jacob}
\end{equation}
Note that $\mbn{X}_{,\xi^\alpha}$ is in fact equivalent to $\mbn{G}_\alpha^0$, however, we try to use the original notation. The derivatives of shape functions with respect to the local nodal coordinate system~$\mbn{A}_I$  are
\begin{equation}
\label{eq:shapefundertransf}
	\begin{bmatrix}
		N_{I,1}\\
		N_{I,2}
	\end{bmatrix}=
		\mbn{J}^{-1}
	\begin{bmatrix}
		N_{I,\xi^1}\\
		N_{I,\xi^2}
	 \end{bmatrix}\,,
\end{equation}
where the left side shows the derivatives of the basis function $N_I$ with respect to the axis~$\alpha$, which is defined by the base vector $\mbn{A}_{\alpha I}$; c.f.~Fig.~\ref{fig_SEMI_element}. We can rearrange the shell strains in Eq.~\eqref{shell_strain_2} in a vector format as 
\begin{equation}\label{shell_strain_vector}
	\mbn{\epsilon}=[\epsilon_{11}, \epsilon_{22}, 2\epsilon_{12}, \kappa_{11}, \kappa_{22},2\kappa_{12},\gamma_1, \gamma_2]\,.
\end{equation} 
Considering Eq.~\eqref{shell_strain_1} and plane stress assumptions, one can calculate the resultants of the second Piola-Kirchhoff stresses throughout the thickness of the shell. These quantities, which are work conjugate to $\mbn{\epsilon}$, are called shell stresses. They are assembled according to \cite{wagner_2005} by
\begin{equation}\label{shell_stress_vector}
	\mbn{\sigma}=[n^{11}, n^{22}, n^{12}, m^{11}, m^{22}, m^{12}, q^1, q^2 ]\,.
\end{equation}
The constitutive relation between shell strains and stresses is
\begin{equation}
	\mbn{\sigma}=\mbn{C} \cdot \mbn{\epsilon}\,,
\end{equation}
where
\begin{equation}
	\mbn{C}=
	\begin{bmatrix}
		h \mbn{C}_p &\mbn{0} &\mbn{0} \\
		\mbn{0} &\frac{h^3}{12} \mbn{C}_p &\mbn{0}\\
		\mbn{0} &\mbn{0} & h \mbn{C}_s 
	\end{bmatrix}
\end{equation}
and
\begin{subequations}
	\begin{align}
	\mbn{C}_p&=\frac{E_e}{1-\nu^2}
	\begin{bmatrix}
		1 &\nu &0 \\
		\nu &1 &0 \\
		0 &0 &\frac{1-\nu}{2}
	\end{bmatrix}\, ,\\[8pt]
	\mbn{C}_s&= k_s \frac{E_e}{2(1+\nu)}
	\begin{bmatrix}
		1 &0\\
		0 &1
	\end{bmatrix}\,.
	\end{align}
\end{subequations}
Here $E_e$ and $\nu$ are Young's modulus and Poisson's ratio, respectively. $k_s$ is the Timoshenko shear modification factor and is taken as 5/6. Using Eqs.~\eqref{shell_strain_vector} and \eqref{shell_stress_vector}, the principle of virtual energy for the shell surface depicted in Fig.~\ref{fig_shell_conf} can be written as
\begin{equation}\label{virtual_work}
	\begin{split}
	G(\mbn{v}, \delta\mbn{v})=&\int_{\Omega_0} \delta\mbn{\epsilon}^{\mathrm{T}} \mbn{\sigma} \, \mathrm{d}A - \int_{\Omega_0} \delta\mbn{v}^{\mathrm{T}} \overline{\mbn{P}} \, \mathrm{d}A - \int_{\Gamma_0} \delta\mbn{v}^{\mathrm{T}} \overline{\mbn{t}}  \, \mathrm{d}s= 0 \, ,
	\end{split}
\end{equation}
where $\overline{\mbn{P}}$ and $\overline{\mbn{t}}$  are conservative surface load and boundary tractions on the surface $\Omega_0$ and the boundary $\Gamma_0$, respectively. Please note that the entries of $\overline{\mbn{P}}$ and $\overline{\mbn{t}}$ are the conjugate forces and moments to the displacements and rotations in $\mbn{v}$, respectively. Eq.~\eqref{virtual_work} is highly nonlinear with respect to the DOFs, i.e.~$\mbn{v}$. Therefore, the Newton-Raphson method is applied to solve the equations, and the linearized form is given by
\begin{equation}\label{newton_1}
	\mathrm{L}(G(\mbn{v}, \delta\mbn{v}))=G + \mathrm{D}_\mathrm{d}(G) \cdot \Delta\mbn{v}\,. 
\end{equation}
The operator $\mathrm{D}_\mathrm{d}$ represent derivatives with respect to $\mbn{v}$, hence it yields
\begin{equation}\label{newton_2}
	\mathrm{D}_\mathrm{d}(G) \cdot \Delta\mbn{v}=\int_{\Omega_0} \delta\mbn{\epsilon}^{\mathrm{T}} \, \mbn{C}\, \Delta{\mbn{\epsilon}} \, \mathrm{d}A+\int_{\Omega_0} \Delta\delta\mbn{\epsilon}^{\mathrm{T}} \, \mbn{\sigma} \,\mathrm{d}A 
\end{equation}
with
\begin{equation}
	\delta(.)=\frac{\partial(.)}{\partial\mbn{v}}\cdot \delta\mbn{v} \quad \mathrm{and} \quad \Delta(.)=\frac{\partial(.)}{\partial\mbn{v}}\cdot \Delta\mbn{v}\,.
\end{equation}
Thus, it is necessary to calculate the first and second derivatives of shell strains and then factorize $\delta\mbn{v}$ and $\Delta\mbn{v}$ to eliminate the former and solve for the latter. Utilizing Eq.~\eqref{shell_strain_2}, the first variations of the interpolated shell strains are
\begin{equation}\label{var_shell_strains_1}
  \delta\mbn{\epsilon}^h=
  \begin{bmatrix}
  	\delta\epsilon_{11}^h \\[0.2cm]
  	\delta\epsilon_{22}^h\\[0.2cm]
  	2 \delta\epsilon_{12}^h \\[0.2cm]
  	\delta\kappa_{11}^h \\[0.2cm]
  	\delta\kappa_{22}^h \\[0.2cm]
  	2 \delta\kappa_{12}^h \\[0.2cm]
  	\gamma_1^h \\[0.2cm]
  	\gamma_2^h
  \end{bmatrix}=
  \begin{bmatrix}
  	\delta\mbn{x}_{,1}^h \cdot \mbn{x}_{,1}^h\\[0.2cm]
  	\delta\mbn{x}_{,2}^h \cdot \mbn{x}_{,2}^h\\[0.2cm]
  	\delta\mbn{x}_{,1}^h \cdot \mbn{x}_{,2}^h + \mbn{x}_{,1}^h \cdot \delta\mbn{x}_{,2}^h\\[0.2cm]
  	\delta\mbn{x}_{,1}^h \cdot \mbn{d}_{,1}^h + \mbn{x}_{,1}^h \cdot \delta\mbn{d}_{,1}^h\\[0.2cm]
  	\delta\mbn{x}_{,2}^h \cdot \mbn{d}_{,2}^h + \mbn{x}_{,2}^h \cdot \delta\mbn{d}_{,2}^h\\[0.2cm]
  	\delta\mbn{x}_{,1}^h \cdot \mbn{d}_{,2}^h + \delta\mbn{x}_{,2}^h \cdot \mbn{d}_{,1}^h + \mbn{x}_{,1}^h \cdot \delta\mbn{d}_{,2}^h + \mbn{x}_{,2}^h \cdot \delta\mbn{d}_{,1}^h\\[0.2cm]
  	\delta\mbn{x}_{,1}^h \cdot \mbn{d}^h + \mbn{x}_{,1}^h \cdot \delta\mbn{d}^h\\[0.2cm]
  	\delta\mbn{x}_{,2}^h \cdot \mbn{d}^h + \mbn{x}_{,2}^h \cdot \delta\mbn{d}^h
  \end{bmatrix}\,.
\end{equation}
The interpolated variation of directors is
\begin{equation} 
		\delta\mbn{d}^h =\sum_{I=1}^{n_{en}} N_I \, \delta\mbn{d}_I \, 
\end{equation}
and the interpolated variation of derivatives of mid-surface position vector and the directors are
\begin{subequations}\label{interpolation_discrete}
	\begin{align}
		&\delta\mbn{x}^h_{,\alpha} =\sum_{I=1}^{n_{en}} N_{I,\alpha} \, \delta\mbn{x}_I \, ,\\
		&\delta\mbn{d}^h_{,\alpha} =\sum_{I=1}^{n_{en}} N_{I,\alpha} \, \delta\mbn{d}_I \, .
	\end{align}
\end{subequations}
Here it is to be noted, that the transformed derivatives of the shape functions, given in Eq.~\eqref{eq:shapefundertransf}, have to be used.
While the deformation parts of the formulation are quite simple, the rotational parts are more complicated. We follow the formulation of \cite{wagner_2005,grutt_2000} to establish the relation between the variation of the directors at the nodes and the variation of the rotational DOFs. This yields 
\begin{equation}
	\delta\mbn{d}_I=\delta\mbn{R}_I \mbn{D}_I=\mbn{W}_I^{\mathrm{T}} \delta\mbn{w}_I \,,
\end{equation}
where
\begin{equation}
	\mbn{W}_I=\mathrm{skew}(\mbn{d}_I)
\end{equation}
and
\begin{equation}\label{wi}
	\delta\mbn{w}_I=\mbn{H}_I \delta\mbn{\omega}_I\,.
\end{equation}
$\mbn{H}_I$ is
\begin{equation}\label{eq_H}
\mbn{H}_I=\mbn{1}+\frac{1-\mathrm{cos}\, \omega_I}{\omega_I^2} \mbn{\Omega}_I+\frac{\omega_I-\mathrm{sin}\, \omega_I}{\omega_I^3} \mbn{\Omega}_I^2 \,.
\end{equation}
Summing up, the variation of the nodal director vector is related by
\begin{equation}
	\delta\mbn{d}_I=\mbn{W}_I^{\mathrm{T}} \mbn{H}_I \delta\mbn{\omega}_I
\end{equation}
to the variation of the axial vector of the rotation $\delta\mbn{\omega}_I$.
Using Eq.~\eqref{beta_to_omega} further introduces the relation
\begin{equation}\label{var_dir_beta}
	\delta\mbn{d}_I=\mbn{T}_I \, \delta\mbn{\beta}_I \qquad \textrm{with} \qquad \mbn{T}_I=\mbn{W}_I^{\mathrm{T}} \mbn{H}_I \mbn{T}_{3I} 
\end{equation}
to the variation of the nodal rotations.
Considering Eqs.~\eqref{var_shell_strains_1}-\eqref{var_dir_beta}, the finite element approximation of variations of the shell strains for one element is 
\begin{equation} \label{B_matrix}
	\begin{split}
	\begin{bmatrix}
		\delta\epsilon_{11}^h \\[0.2cm]
		\delta\epsilon_{22}^h\\[0.2cm]
		2 \delta\epsilon_{12}^h \\[0.2cm]
		\delta\kappa_{11}^h \\[0.2cm]
		\delta\kappa_{22}^h \\[0.2cm]
		2 \delta\kappa_{12}^h \\[0.2cm]
		\gamma_1^h \\[0.2cm]
		\gamma_2^h
	\end{bmatrix}  &= \sum_{I=1}^{n_{en}}
	\begin{bmatrix}
		N_{I,1} {\mbn{x}_{,1}^h}^\mathrm{T} &\mbn{0}\\[0.2cm]
		N_{I,2} {\mbn{x}_{,2}^h}^\mathrm{T} &\mbn{0}\\[0.2cm]
		N_{I,1} {\mbn{x}_{,2}^h}^\mathrm{T}+N_{I,2} {\mbn{x}_{,1}^h}^\mathrm{T} &\mbn{0}\\[0.2cm]
		N_{I,1} {\mbn{d}_{,1}^h}^\mathrm{T} & N_{I,1} ({\mbn{x}_{,1}^h}^{\mathrm{T}} \, \mbn{T}_I)\\[0.2cm]
	    N_{I,2} {\mbn{d}_{,2}^h}^{\mathrm{T}} & N_{I,2} ({\mbn{x}_{,2}^h}^\mathrm{T} \, \mbn{T}_I)\\[0.2cm]
	    N_{I,1} {\mbn{d}_{,2}^h}^\mathrm{T}+N_{I,2} {\mbn{d}_{,1}^h}^\mathrm{T} & N_{I,1}{(\mbn{x}_{,2}^h}^\mathrm{T} \, \mbn{T}_I)+N_{I,2}({\mbn{x}_{,1}^h}^{\mathrm{T}} \, \mbn{T}_I)\\[0.2cm]
	    N_{I,1} {\mbn{d}^h}^\mathrm{T} & N_I ({\mbn{x}_{,1}^h}^\mathrm{T} \, \mbn{T}_I)\\[0.2cm]
	    N_{I,2} {\mbn{d}^h}^{\mathrm{T}} & N_I({\mbn{x}_{,2}^h}^{\mathrm{T}} \,\mbn{T}_I)
	\end{bmatrix}
	\begin{bmatrix}
		\delta\mbn{u}_I \\
		\delta\mbn{\beta}_I
	\end{bmatrix}  \\[0.3cm]
        & = \sum_{I=1}^{n_{en}} \mbn{B}_I \delta\mbn{v}_I
	= \mbn{B}^e \delta\mbn{v}^e\,,
	\end{split}
\end{equation}
where $\mbn{v}^e$ is the vector of DOFs for one element. The second variations of shell strains are 
\begin{equation}\label{sec_var_shell_strains_1}
	\Delta\delta\mbn{\epsilon}^h=
	\begin{bmatrix}
		\Delta\delta\epsilon_{11}^h \\[0.2cm]
		\Delta\delta\epsilon_{22}^h\\[0.2cm]
		2 \Delta\delta\epsilon_{12}^h \\[0.2cm]
		\Delta\delta\kappa_{11}^h \\[0.2cm]
		\Delta\delta\kappa_{22}^h \\[0.2cm]
		2 \Delta\delta\kappa_{12}^h \\[0.2cm]
		\Delta\gamma_1^h \\[0.2cm]
		\Delta\gamma_2^h
	\end{bmatrix}=
	\begin{bmatrix}
		\delta\mbn{x}_{,1}^h \cdot \Delta\mbn{x}_{,1}^h\\[0.2cm]
		\delta\mbn{x}_{,2}^h \cdot \Delta\mbn{x}_{,2}^h\\[0.2cm]
		\delta\mbn{x}_{,1}^h \cdot \Delta\mbn{x}_{,2}^h + \Delta\mbn{x}_{,1}^h \cdot \delta\mbn{x}_{,2}^h\\[0.2cm]
		\delta\mbn{x}_{,1}^h \cdot \Delta\mbn{d}_{,1}^h + \Delta\mbn{x}_{,1}^h \cdot \delta\mbn{d}_{,1}^h + \mbn{x}_{,1}^h \cdot \Delta\delta\mbn{d}_{,1}^h\\[0.2cm]
		\delta\mbn{x}_{,2}^h \cdot \Delta\mbn{d}_{,2}^h + \Delta\mbn{x}_{,2}^h \cdot \delta\mbn{d}_{,2}^h + \mbn{x}_{,2}^h \cdot \Delta\delta\mbn{d}_{,2}^h\\[0.2cm]
		\delta\mbn{x}_{,1}^h \cdot \Delta\mbn{d}_{,2}^h + \delta\mbn{x}_{,2}^h \cdot \Delta\mbn{d}_{,1}^h + \Delta\mbn{x}_{,1}^h \cdot \delta\mbn{d}_{,2}^h + \\ 
		\Delta\mbn{x}_{,2}^h \cdot \delta\mbn{d}_{,1}^h + 
		 \mbn{x}_{,1}^h \cdot \Delta\delta\mbn{d}_{,2}^h + \mbn{x}_{,2}^h \cdot \Delta\delta\mbn{d}_{,1}^h\\[0.2cm]
		\delta\mbn{x}_{,1}^h \cdot \Delta\mbn{d}^h + \Delta\mbn{x}_{,1}^h \cdot \delta\mbn{d}^h + \mbn{x}_{,1}^h \cdot \Delta\delta\mbn{d}^h\\[0.2cm]
		\delta\mbn{x}_{,2}^h \cdot \Delta\mbn{d}^h + \Delta\mbn{x}_{,2}^h \cdot \delta\mbn{d}^h + \mbn{x}_{,2}^h \cdot \Delta\delta\mbn{d}^h
	\end{bmatrix}
\end{equation}
In the calculation of Eq.~\eqref{sec_var_shell_strains_1}, we consider that $\delta\mbn{x}_{,\alpha}^h $ is not a function of DOFs, while as can be deduced from Eq.~\eqref{var_dir_beta}, $\delta\mbn{d}^h$ and its derivatives are nonlinear functions of the rotational DOFs. The finite element approximation of the second variation of director and its derivatives are
\begin{equation} \label{sec_var_dir}
	\Delta\delta\mbn{d}^h =\sum_{I=1}^{n_{en}} N_I \, \Delta\delta\mbn{d}_I \, 
\end{equation}
and
\begin{equation}\label{sec_var_der_dir}
		\Delta\delta\mbn{d}^h_{,\alpha} =\sum_{I=1}^{n_{en}} N_{I,\alpha} \, \Delta\delta\mbn{d}_I \, .
\end{equation}
Considering Eqs.~\eqref{sec_var_shell_strains_1}, \eqref{sec_var_dir} and \eqref{sec_var_der_dir}, we need to calculate $\mbn{h} \cdot \Delta\delta\mbn{d}_I$ where $\mbn{h}$ is an arbitrary vector and the representative of the derivatives of $\mbn{x}_{,\alpha}$. From \cite{wagner_2005}, we have
\begin{equation}
	\mbn{h} \cdot \Delta\delta\mbn{d}_I = \delta\mbn{w}_I^{\mathrm{T}} \mbn{M}_I(\mbn{h}) \Delta\mbn{w}_I\,,
\end{equation}
where $\delta\mbn{w}_I$ is defined in Eq.~\eqref{wi} and $\mbn{M}_I(\mbn{h})$ is computed by
\begin{equation}\label{eq_MI}
	\begin{split}
	\mbn{M}_I(\mbn{h}) &= \frac{1}{2}(\mbn{d}_I \otimes \mbn{h} + \mbn{h} \otimes \mbn{d}_I) + \frac{1}{2}(\mbn{t}_I \otimes \mbn{\omega}_I + \mbn{\omega}_I \otimes \mbn{t}_I) + c_{10} \mbn{1}\\[4pt]
	\mbn{t}_I &= -c_3 \mbn{b}_I + c_{11} (\mbn{b}_I \cdot \mbn{\omega}_I)\mbn{\omega}_I\\[4pt]
	\mbn{b}_I &= \mbn{d}_I \times \mbn{h}\\[4pt]
	c_3 &=\frac{\omega_I \,\mathrm{sin}\, \omega_I + 2(\mathrm{cos}\,\omega_I - 1)}{\omega_I^2(\mathrm{cos}\,\omega_I - 1)}\\[4pt]
	c_{10}&=\bar{c}_{10}(\mbn{b}_I \cdot \mbn{\omega}_I)-(\mbn{d}_I \cdot \mbn{h})\\[3pt]
	\bar{c}_{10} &= \frac{\mathrm{sin}\,\omega_I - \omega_I}{2 \, \omega_I (\mathrm{cos}\,\omega_I - 1)}\\[4pt]
	c_{11} &= \frac{4(\mathrm{cos}\,\omega_I-1) + \omega_I^2 + \omega_I \, \mathrm{sin}\,\omega_I}{2\, \omega_I^4(\mathrm{cos}\,\omega_I-1)}\,.
	\end{split}
\end{equation}
If $\omega_I$ approaches zero, the calculation of the coefficients may contain some numerical singularities. Therefore, the Taylor series expansion of them should be used instead of their closed form expression \cite{wagner_2005}. Regarding Eq.~\eqref{sec_var_shell_strains_1}, the finite element approximation of the second variation of the shell strains are
\begin{equation}\label{secvar_2}
	\begin{split}
    \Delta\delta\epsilon_{\alpha\beta}^h =& \frac{1}{2} \sum_{I=1}^{n_{en}} \sum_{K=1}^{n_{en}} \delta\mbn{u}_I ^{\mathrm{T}}(N_{I,\alpha} N_{K,\beta}+N_{I,\beta}N_{K,\alpha}) \mbn{1}_{3\times3} \Delta\mbn{u}_K \\
	\Delta\delta\kappa_{\alpha \beta}^h =& \frac{1}{2}  \sum_{I=1}^{n_{en}} \sum_{K=1}^{n_{en}} \Big\{\delta\mbn{u}_I^{\mathrm{T}} (N_{I,\alpha} N_{K,\beta} + N_{I,\beta} N_{K,\alpha}) \mbn{T}_K \Delta\mbn{\beta}_K \\
	&+\delta\mbn{\beta}_I^{\mathrm{T}} \mbn{T}_I^{\mathrm{T}} (N_{I,\alpha} N_{K,\beta}+N_{I,\beta} N_{K,\alpha}) \Delta\mbn{u}_K \\
	&+ \delta_{IK} \, \delta\mbn{\beta}_I^{\mathrm{T}} \mbn{T}_{3I}^{\mathrm{T}} \mbn{H}_I^{\mathrm{T}} \left[\mbn{M}_I(N_{I,\alpha}  \, \mbn{x}_{,\beta})+\mbn{M}_I(N_{I,\beta} \, \mbn{x}_{,\alpha})\right] \mbn{H}_K \, \mbn{T}_{3K} \Delta\mbn{\beta}_K \Big\}\\
	\Delta\delta\gamma_\alpha^h=&\sum_{I=1}^{n_{en}} \sum_{K=1}^{n_{en}} \delta\mbn{u}_I^{\mathrm{T}} N_{I,\alpha} N_K \mbn{T}_K \Delta\mbn{\beta}_K+\delta\mbn{\beta}_I^\mathrm{T} \mbn{T}_I^\mathrm{T} N_I N_{K,\alpha} \Delta\mbn{u}_K\\
	&+ \delta_{IK} \, \delta\mbn{\beta}_I^\mathrm{T} \mbn{T}_{3I}^\mathrm{T} \mbn{H}_I^\mathrm{T}  \,\mbn{M}_I(N_I \, \mbn{x}_{,\alpha}) \mbn{H}_K \, \mbn{T}_{3K} \Delta{\mbn{\beta}}_K\,,
	\end{split}
\end{equation} 
where $\delta_{IK}$ is the Kronecker delta; see \cite{Dornisch.2015b,wagner_2005}. Finally, by utilizing Eq.~\eqref{secvar_2}, the term $\Delta\delta\mbn{\epsilon}^{\mathrm{T}} \mbn{\sigma}$ of Eq.~\eqref{newton_2} can be written as follows
\begin{equation}
	\begin{split}
		{\Delta\delta\mbn{\epsilon}^h}^{\mathrm{T}} \mbn{\sigma}&=\sum_{I=1}^{n_{en}}\sum_{K=1}^{n_{en}} \delta\mbn{v}_I {\mbn{k}^G_{IK}} \Delta\mbn{v}_K=\\
		&\sum_{I=1}^{n_{en}}\sum_{K=1}^{n_{en}}\\
		&\begin{bmatrix}
			\delta\mbn{u}_I\\
			\delta\mbn{\beta}_I
		\end{bmatrix}^\mathrm{T}
		\begin{bmatrix}
			\hat{n}_{IK}\mbn{1}_{3\times3}
			&(\hat{m}_{IK}+\hat{q}_{IK}) \mbn{T}_K\\[0.2cm]
		     \mbn{T}_I^\mathrm{T} (\hat{m}_{IK}+\widetilde{q}_{IK})
		     &\delta_{IK} \, \mbn{T}_{3I}^\mathrm{T} \mbn{H}_I^\mathrm{T}  \hat{\mbn{M}}_I \,  \mbn{H}_K \, \mbn{T}_{3K}
		\end{bmatrix}
		\begin{bmatrix}
			 \Delta\mbn{u}_K\\
			\Delta\mbn{\beta}_K
		\end{bmatrix}\,,
	\end{split}
    \label{eq:k_g_ik}
\end{equation}
where
\begin{equation}
\label{eq:kgik_entries}
	\begin{split}
		\hat{n}_{IK} &= n^{11} N_{I,1} N_{K,1}+n^{22} N_{I,2} N_{K,2}+n^{12}(N_{I,1} N_{K,2}+N_{I,2}N_{K,1})\\
		\hat{m}_{IK} &= m^{11} N_{I,1} N_{K,1}+m^{22} N_{I,2} N_{K,2}+m^{12}(N_{I,1} N_{K,2}+N_{I,2}N_{K,1})\\
		\hat{q}_{IK} &= q^1 N_{I,1} N_K+q^2  N_{I,2} N_K\\
		\widetilde{q}_{IK} &= q^1 N_I N_{K,1}+q^2 N_I N_{K,2}\\
        \hat{\mbn{h}}_I &= \sum_{\alpha=1}^{2}{\left(m^{1\alpha}N_{I,1}+m^{\alpha2} N_{I,2}+q^\alpha N_I\right) \,\mbn{x}^h_{,\alpha}}\\
        \hat{\mbn{M}}_I&=\mbn{M}_I(\hat{\mbn{h}}_I)\,.
	\end{split}
\end{equation}
\begin{remark}
   It is to be mentioned that the computation of the lower right part of $ {\mbn{k}^G_{IK}}$ constitutes a significant part of the total matrix assembly time.
   Since $\mbn{M}_I(\mbn{h})$ is linear in $\mbn{h}$, there are several possibilities to compute~$\hat{\mbn{M}}_I$. The number of required multiplications is minimal if the version presented here is chosen. See also a discussion on a similar case in \cite[pp. 102-103]{Dornisch.2015b}.   
    
    
    
    
\end{remark}
By using $ {\mbn{k}^G_{IK}}$, which corresponds to node $I$ and $K$ and taking the integration, the geometrical stiffness matrix $\mbn{k}_G$ for the entire element can be easily constructed accordingly. Consequently, for an element with surface $\Omega_0^e$ and boundary $\Gamma_0^e$, the element stiffness matrix, $\mbn{k}^t$, element load vectors $\mbn{f}_s^e$ and $\mbn{f}_l^e$ are
\begin{equation}\label{basic-matrices}
	\begin{split}
         \mbn{k}_E=&\int_{\Omega^e_0} {\mbn{B}^e}^{\mathrm{T}} \mbn{C} \, \mbn{B}^e \mathrm{d}A\\
         \mbn{k}^t=&\mbn{k}_E+\mbn{k}_G\\
	\mbn{f}_s^e=&\int_{\Omega_0^e}\mbn{N}^{\mathrm{T}} \, \overline{\mbn{P}} \, \mathrm{d}A \, ,\\
	\mbn{f}_l^e=& \int_{\Gamma_0^e} \mbn{N}^{\mathrm{T}} \, \overline{\mbn{t}} \, \mathrm{d}s \,,
	\end{split}
\end{equation}
where the basis function matrix
\begin{equation} \label{N_matrix}
	\mbn{N} =
	\begin{bmatrix}
    N_1\mbn{1}_{5\times5} & N_2\mbn{1}_{5\times5} & \ldots & N_{n_{en}}\mbn{1}_{5\times5}
	\end{bmatrix}
\end{equation}
is used. For elements at the shell intersections, Eq.~\eqref{N_matrix} must be modified accordingly to account for the drilling rotation DOFs. The surface integrals in Eq.~\eqref{basic-matrices} are computed numerically by using the GLL method and
\begin{equation}\label{area_ratio_1}
	\mathrm{d}A=|\mbn{X}_{,\xi^1}^h \times \mbn{X}_{,\xi^2}^h| \, \mathrm{d}\xi^1 \mathrm{d}\xi^2 \,.
\end{equation}
Regarding Eq.~\eqref{newton_1}, $G$ is calculated by using Eq.~\eqref{B_matrix} and finally the residual force is
\begin{equation}
\label{eq:res_final}
	\mbn{f}^e_{res} = \mbn{f}_s^e+\mbn{f}_l^e - \int_{\Omega_0^e} {\mbn{B}^e}^T \mbn{\sigma} \,\mathrm{d}A\,.
\end{equation} 
Global stiffness matrix and residual force are assembled from their local counterparts provided in Eqs.~\eqref{basic-matrices} and~\eqref{eq:res_final} using standard finite element routines. Afterwards, they are used in the Newton-Raphson procedure in accordance with Eqs.~\eqref{virtual_work} and \eqref{newton_1}.

\subsubsection{Spectral Element Method, Non-Isoparametric formulation (SEMN)}\label{section_SEMN}
\begin{figure}[ht] 
	\center{\includegraphics[scale=0.5]{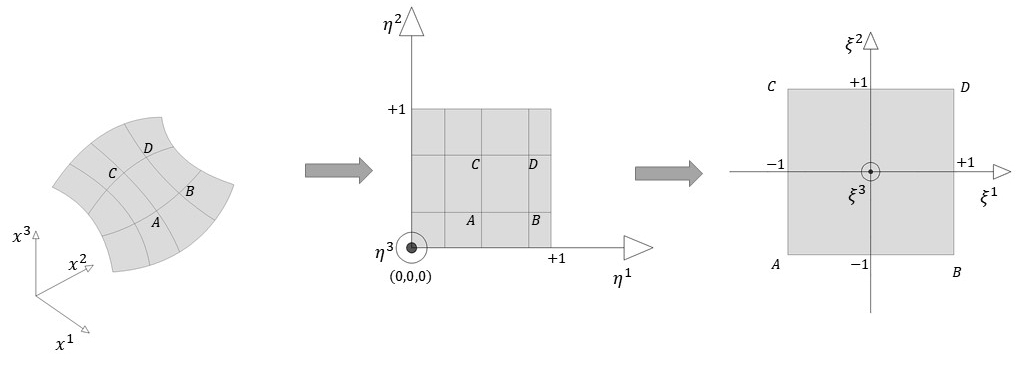} }
	\caption{\centering Physical, first (IGA) and second (spectral element) parametric spaces. Depicted mesh is according to the first four scenarios.}\label{fig_maps_semn}
\end{figure}

In SEMN, as in SEMI, Lagrange shape functions with Lobatto nodes are used to interpolate the unknown fields, while control point coordinates, orders and knot vectors are directly imported from Computer-Aided design (CAD) files to define the geometry \cite{nima_2024}. Therefore, SEMN benefits from the exact geometry definition, similar to isogeometric analysis (IGA) \cite{hughes_2005}. It means that we attach a spectral mesh to the NURBS based geometry. Meshing can be based on five scenarios: 
\begin{enumerate}
    \item elements coincide with the non-zero-length elements of IGA
    \item internal borders of elements coincide only at the internal knots with predefined number of repetition (according to the continuity of the shell) 
    \item element borders at any arbitrary knot 
    \item a mixture of all three previous cases 
    \item Unstructured meshing, it means that the element borders are not necessarily parallel with the $\eta^1$ and $\eta^2$ axes anymore. 
    \end{enumerate} 
    So, instead of the single mapping from parametric space to physical space for SEMI presented in Fig.~\ref{fig_map_semi}, we will use two mappings for SEMN; see Fig.~\ref{fig_maps_semn}.
The parametric space, with coordinates $\eta^1$ and $\eta^2$, is used to define NURBS shape functions. The second space, with coordinates $\xi^1$ and $\xi^2$, is used to compute the Lagrange shape functions with Lobatto nodes. To introduce the notation, a brief overview of the NURBS shape functions is presented here. Univariate B-spline basis functions are computed recursively by the Cox-de Boor formula
\begin{align}\label{bspline_def_1}
	p=0: &\bar{N}^0_i=
	\begin{cases}
		1 & \text{if  \( \eta_i\leq \eta \leq \eta_{i+1} \) } \\
		0 & \text{otherwise}
	\end{cases}\\
    p>0: &\bar{N}^p_i(\eta)= \frac{\eta-\eta_i}{\eta_{i+p} - \eta_i}\bar{N}^{p-1}_i(\eta)+\frac{\eta_{i+p+1}-\eta}{\eta_{i+p+1}-\eta_{i+1}} \bar{N}^{p-1}_{i+1}(\eta)
\end{align}
based on the open and non-decreasing knot vector $\mbn{\eta}=\{\eta_1, \eta_2,...,\eta_{n+p+1}\}$, where $p$ is the order of the B-spline basis function, $n$ is the number of basis functions and $\eta\in[\eta_1,\eta_{n+p+1}]$. For more detail, see \cite{Cottrell.2009,Piegl.1997}. The coordinate of a point, $\mbn{X}_\textrm{NURBS}$, on the NURBS mid-surface of a shell patch can be defined by utilizing a net of $n^1\times n^2$ control points, the product of univariate B-spline basis functions of orders $p^1$ in $\eta^1$ and $p^2$ in $\eta^2$ directions and weight coefficients $w_{ij}$. This results in
\begin{equation}
	\mbn{X}_{\mathrm{NURBS}}(\eta^1, \eta^2) = \sum_{L= 1}^{n_{en}}\hat{N}_L(\eta^1,\eta^2) \mbn{X}_L
\end{equation}
with 
\begin{equation}
	\hat{N}_L=\frac{\bar{N}^{p^1}_i(\eta^1)\bar{N}^{p^2}_j(\eta^2)w_{ij}}{W(\eta^1,\eta^2)}
\end{equation}
and
\begin{equation}
	W(\eta^1, \eta^2)= \sum_{i=1}^{n^1}{} \sum_{j=1}^{n^2} \bar{N}^{p^1}_i(\eta^1)\bar{N}^{p^2}_j(\eta^2)w_{ij}.
\end{equation}
 The number of 2D shape functions that have impact on the arbitrary point with the coordinate $\eta^1$ and $\eta^2$  in the parametric space is denoted by $n_{en} = (1+p^1)(1+p^2)$. $\mbn{X}_L$ is the coordinate vector of the \textit{L}-th control point in the set of the control points which has impact on the arbitrary point. To be used in the next parts, the detailed algorithm to calculate the derivatives of NURBS functions can be found in \cite{Piegl.1997}.
\par The concept of SEMN shows that most parts of the formulation in Sec. \ref{semi} are valid here, requiring only a slight modification for the calculation of the Jacobian matrix. The Jacobian matrix of the whole mapping from the physical space into the second parametric space is comprised of two parts, $\mbn{J}_1$ from physical to IGA parametric space and $\mbn{J}_2$ from IGA to spectral parametric space. Following the same logic, Eq.~\eqref{eq:jacob} can be adjusted to
\begin{equation}
	\mbn{J}_1=
	\begin{bmatrix} \label{jacobi_1_semn}
		\mbn{X}_{\mathrm{NURBS},\,\eta^1} \cdot \mbn{A}_1 & \mbn{X}_{\mathrm{NURBS},\,\eta^1} \cdot \mbn{A}_2\\
		\mbn{X}_{\mathrm{NURBS},\,\eta^2} \cdot \mbn{A}_1 & \mbn{X}_{\mathrm{NURBS}        ,\,\eta^2} \cdot \mbn{A}_2
	\end{bmatrix}.
\end{equation}
The Jacobian matrix of the transformation from IGA to SEM parametric spaces can be written as
\begin{equation} \label{j2}
	(\boldsymbol{\mathrm{J}}_2)_{ij}= \left[\frac{\partial \eta^j}{\partial \xi^i}\right].
\end{equation}
This is based on simple linear Lagrange shape functions used to map element geometries from the IGA to the spectral element parametric space. Finally, the Jacobian of the entire transformation from physical to SEM parametric space is
\begin{equation}\label{entirejac}
	\mbn{J}_{\mathrm{total}}= \boldsymbol{\mathrm{J}}_2 \;\boldsymbol{\mathrm{J}}_1.
\end{equation}
and for SEMN
\begin{equation}
	\begin{bmatrix}
		N_{I,1}\\
		N_{I,2}
	\end{bmatrix}=
	\mbn{J}_{\mathrm{total}}^{-1}
	\begin{bmatrix}
		N_{I,\xi^1}\\
		N_{I,\xi^2}
	\end{bmatrix}.
\end{equation}
has to be used instead of Eq.~\eqref{eq:shapefundertransf}.
In the numerical integration, the GLL quadrature is used and instead of Eq.~\eqref{area_ratio_1}, the following equation is utilized
\begin{equation}
	\mathrm{d}A=|\mbn{X}_{\mathrm{NURBS},\eta^1}^h \times \mbn{X}_{\mathrm{NURBS},\eta^2}^h| \,
   \left|\mbn{J}_2\right| \mathrm{d}\xi^1 \mathrm{d}\xi^2 \,.
\end{equation}
More details about SEMN can be found in Ref.~\cite{nima_2024}, where a linear Reissner--Mindlin SEMN shell formulation was proposed.

\subsubsection{Isogeometric Analysis Method (IGA)}
In IGA, the same shape functions, typically NURBS, used to generate the geometry are also utilized to interpolate unknown fields \cite{hughes_2005}. We refrain from delving further into details about IGA and focus solely on the definitions and formulations relevant to this contribution.
\par In IGA, the geometry is represented exactly as designed in a CAD tool. Thus, we can directly compute perfectly aligned local basis systems at the integration points. However, the shape of the geometry is governed by control points, which in general do not lie on the shell surfaces and can lie far away from elements upon which they have influence. Since nodal DOFs are defined at control points, rotation-based isogeometric Reissner--Mindlin shell formulations require nodal basis system at the control points. As these points do not lie on the physical shell surface, there is no unique and exact choice for control point basis systems; see \cite{Dornisch.2015b}.  Therefore, firstly, a method to calculate local coordinate systems at control points based on the geometry of the surface needs to be implemented. Thereby, the interpolated basis vectors in thickness direction, i.e. $\mbn{A}^h_3=\sum\mbn{a}_{3I}$, should be equivalent to the local thickness, expressed by the director vector $\mbn{D}$, in order to prevent artifical thinning in the reference configuration. Secondly, as in IGA shell elements, control points do not coincide with integration points, interpolated directors at integration points may suffer from artificial thinning during deformation, which may affect profoundly the accuracy of the results of rotation-based Reissner--Mindlin shell formulations \cite{Dornisch.2013}. A potential solution for the first point is to use a least-square projection to compute the initial nodal basis systems at the control points as proposed in \cite{Dornisch.2013}. This method minimizes artificial thinning in the reference configuration, which is getting more prominent with rising orders, if a simpler method like the closest point projection approach of \cite{Benson.2010} is used. To prevent negative influence of artificial thinning, we use the least-square projection method for all IGA computations in our numerical examples. A potential solution for the second point is to use continuous interpolation of the rotation of the current director vector as proposed in \cite{dor_2016}. However, this comes with some shortcomings, which are
\begin{enumerate}
    \item Effort to compute the initial nodal basis systems
    \item Significantly higher numerical effort of the continuous director rotation interpolation in comparison to the standard discrete director rotation interpolation
    \item Introduction of slight path-dependency by the continuous director rotation interpolation; see~\cite{dor_2016}
\end{enumerate}
Since the main motivation of using spectral elements in our research is the possibility to use  the conventional and more simple standard discrete director rotation interpolation described in Sec.~\ref{semi}, we will compare our new spectral Reissner--Mindlin shell formulation with an IGA shell using both the continuous and the discrete rotational approach. To clarify notation, we will shortly define both approaches.

\paragraph{Isogeometric Reissner--Mindlin shell, Discrete Rotation (IGA-RMD)}
In the standard rotational formulation the current director vector and its variations and derivatives are computed as in our proposed formulation in Eqs.~\eqref{d_inter} and \eqref{d_inter_der}, with the only difference that NURBS basis functions $\hat{N}_I$ are used for the interpolation. This yields
\begin{equation}\label{d_iga_inter}
		\mbn{d}^h =\sum_{I=1}^{n_{en}} \hat{N}_I \, \mbn{d}_I \qquad\textrm{and}\qquad
		\mbn{d}^h_{,\alpha} =\sum_{I=1}^{n_{en}} \hat{N}_{I,\alpha} \, \mbn{d}_I 
\end{equation}
for the interpolated rotated director and its derivatives, respectively. Numerical studies in~\cite{Dornisch.2015b} have shown that this rotational formulation yields diverging deformation results for order elevation. It should only be employed for low order IGA computations. However, since our new SEMI shell uses the same rotational concept, we provide results of this rotational formulation -- denoted by IGA-RMD -- in the upcoming numerical examples.
\paragraph{Isogeometric Reissner--Mindlin shell, Continuous Rotation (IGA-RMC)}

The basic idea of the continuous rotational formulation, initially proposed in~\cite{Dornisch.2013} in the frame of IGA, is to compute the current director vector by a rotation of the interpolated (or exact) reference director vector
\begin{equation}\label{d_iga_inter_2}
	\mbn{d}^h =\mbn{R} \sum_{I=1}^{n_{en}} \hat{N}_I \, \mbn{D}_I \, . 
\end{equation}
There the rotation matrix $\mbn{R}=\mbn{R}(\mbn{\omega}^h)$ is obtained from the interpolated axial vector of the rotation
\begin{equation}
    \mbn{\omega}^h =\sum_{I=1}^{n_{en}} \hat{N}_I \, \mbn{\omega}_I \, . 
\end{equation}
In order to reduce the path dependency of the formulation, which is a negative impact of using the continuous rotation interpolation, in~\cite{dor_2016} it has been proposed to use a multiplicative update of the rotational state. In the multiplicative version, the current rotational matrix $\mbn{R}^{i}$ is computed by 
\begin{equation}
    \mbn{R}^{i}=\Delta\mbn{R}\mbn{R}^{i-1}\qquad\textrm{with}\qquad \Delta\mbn{R}=\mbn{R}(\Delta\mbn{\omega}^{i})
\end{equation}
from the rotation matrix of the last iteration $\mbn{R}^{i-1}$. The term $\mbn{R}(\Delta\mbn{\omega}^{i})$ states that the interpolated incremental axial vector of the rotation $\Delta\mbn{\omega}^{i}$ has to be inserted into the rotation formula, i.e., Eq.~\eqref{rodrig}.
By using the continuous approach, the derivatives of the current director vector 
\begin{equation} \label{d_iga_inter_der_2}
	\mbn{d}^h_{,\alpha} =\mbn{R}_{,\alpha} \sum_{I=1}^{n_{en}} \hat{N}_{I} \, \mbn{D}_I \,+\mbn{R} \sum_{I=1}^{n_{en}} \hat{N}_{I,\alpha} \, \mbn{D}_I \,,
\end{equation}
and even more the variations thereof, get significantly more lengthy than in the discrete approach.
The full formulation and a thorough overview of different rotational formulations can be found in~\cite{Dornisch.2015b}.  There it is shown, that for IGA the more costly continuous rotation approach with multiplicative update yields superior accuracy and proper convergence behavior for order elevation. Thus, we also compare the results of our new SEMI shell to computations of this IGA shell, coined IGA-RMC in the upcoming numerical examples.

\section{A Brief Discussion on the Effectiveness of SEM methods} \label{sec_3}
\begin{figure}
\begin{center}
	\def\svgwidth{0.4\textwidth}
\begingroup%
  \makeatletter%
  \providecommand\color[2][]{%
    \errmessage{(Inkscape) Color is used for the text in Inkscape, but the package 'color.sty' is not loaded}%
    \renewcommand\color[2][]{}%
  }%
  \providecommand\transparent[1]{%
    \errmessage{(Inkscape) Transparency is used (non-zero) for the text in Inkscape, but the package 'transparent.sty' is not loaded}%
    \renewcommand\transparent[1]{}%
  }%
  \providecommand\rotatebox[2]{#2}%
  \newcommand*\fsize{\dimexpr\f@size pt\relax}%
  \newcommand*\lineheight[1]{\fontsize{\fsize}{#1\fsize}\selectfont}%
  \ifx\svgwidth\undefined%
    \setlength{\unitlength}{242.07818219bp}%
    \ifx\svgscale\undefined%
      \relax%
    \else%
      \setlength{\unitlength}{\unitlength * \real{\svgscale}}%
    \fi%
  \else%
    \setlength{\unitlength}{\svgwidth}%
  \fi%
  \global\let\svgwidth\undefined%
  \global\let\svgscale\undefined%
  \makeatother%
  \begin{picture}(1,1.00150891)%
    \lineheight{1}%
    \setlength\tabcolsep{0pt}%
    \put(0,0){\includegraphics[width=\unitlength]{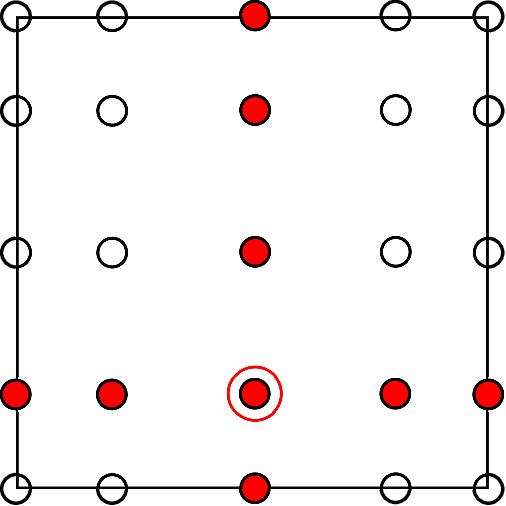}}%
    \put(0.50406528,0.28903783){\color[rgb]{0,0,0}\makebox(0,0)[lt]{\lineheight{1.25}\smash{\begin{tabular}[t]{l}$Q$\end{tabular}}}}%
  \end{picture}%
\endgroup%

\end{center}
\caption{\centering The cross pattern: All nodes of the set $\eta_Q$, i.e. nodes which have non-zero $N_I$ or $N_{I,\alpha}$ at the integration point $Q$, are shown by red color for an element with order $p=4$, as an instance. \label{fig_cross_pattern}}
\end{figure}

\begin{figure}[ht]
	\begin{subfigure}[t]{\linewidth}
            \centering
            \def\svgwidth{0.4\textwidth}
\begingroup%
  \makeatletter%
  \providecommand\color[2][]{%
    \errmessage{(Inkscape) Color is used for the text in Inkscape, but the package 'color.sty' is not loaded}%
    \renewcommand\color[2][]{}%
  }%
  \providecommand\transparent[1]{%
    \errmessage{(Inkscape) Transparency is used (non-zero) for the text in Inkscape, but the package 'transparent.sty' is not loaded}%
    \renewcommand\transparent[1]{}%
  }%
  \providecommand\rotatebox[2]{#2}%
  \newcommand*\fsize{\dimexpr\f@size pt\relax}%
  \newcommand*\lineheight[1]{\fontsize{\fsize}{#1\fsize}\selectfont}%
  \ifx\svgwidth\undefined%
    \setlength{\unitlength}{499.43989563bp}%
    \ifx\svgscale\undefined%
      \relax%
    \else%
      \setlength{\unitlength}{\unitlength * \real{\svgscale}}%
    \fi%
  \else%
    \setlength{\unitlength}{\svgwidth}%
  \fi%
  \global\let\svgwidth\undefined%
  \global\let\svgscale\undefined%
  \makeatother%
  \begin{picture}(1,0.30081699)%
    \lineheight{1}%
    \setlength\tabcolsep{0pt}%
    \put(0,0){\includegraphics[width=\unitlength]{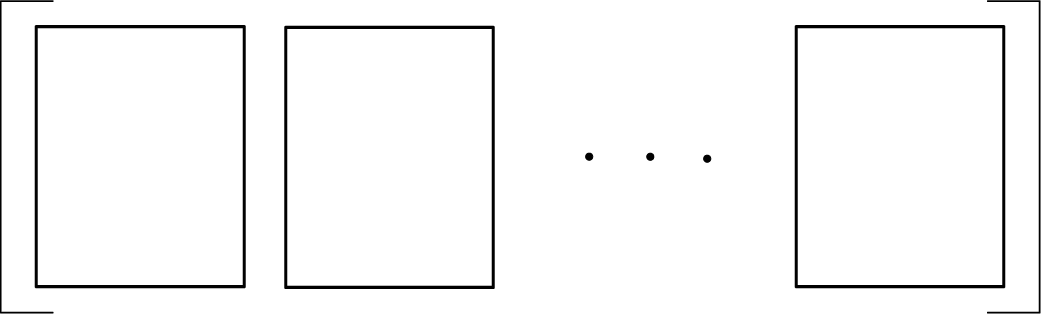}}%
    \put(0.08457919,0.14361355){\color[rgb]{0,0,0}\makebox(0,0)[lt]{\lineheight{1.25}\smash{\begin{tabular}[t]{l}$\mbn{B}_1$\end{tabular}}}}%
    \put(0.32412754,0.14286381){\color[rgb]{0,0,0}\makebox(0,0)[lt]{\lineheight{1.25}\smash{\begin{tabular}[t]{l}$\mbn{B}_2$\end{tabular}}}}%
    \put(0.78863499,0.14361355){\color[rgb]{0,0,0}\makebox(0,0)[lt]{\lineheight{1.25}\smash{\begin{tabular}[t]{l}$\mbn{B}_{n_{en}}$\end{tabular}}}}%
  \end{picture}%
\endgroup%

                \caption{$\mbn{B}_I$ to $\mbn{B}^e$}        
       \end{subfigure}\\

       \begin{subfigure}[t]{\linewidth}
               \centering
               \def\svgwidth{0.4\textwidth}
\begingroup%
  \makeatletter%
  \providecommand\color[2][]{%
    \errmessage{(Inkscape) Color is used for the text in Inkscape, but the package 'color.sty' is not loaded}%
    \renewcommand\color[2][]{}%
  }%
  \providecommand\transparent[1]{%
    \errmessage{(Inkscape) Transparency is used (non-zero) for the text in Inkscape, but the package 'transparent.sty' is not loaded}%
    \renewcommand\transparent[1]{}%
  }%
  \providecommand\rotatebox[2]{#2}%
  \newcommand*\fsize{\dimexpr\f@size pt\relax}%
  \newcommand*\lineheight[1]{\fontsize{\fsize}{#1\fsize}\selectfont}%
  \ifx\svgwidth\undefined%
    \setlength{\unitlength}{490.55996704bp}%
    \ifx\svgscale\undefined%
      \relax%
    \else%
      \setlength{\unitlength}{\unitlength * \real{\svgscale}}%
    \fi%
  \else%
    \setlength{\unitlength}{\svgwidth}%
  \fi%
  \global\let\svgwidth\undefined%
  \global\let\svgscale\undefined%
  \makeatother%
  \begin{picture}(1,1.00048925)%
    \lineheight{1}%
    \setlength\tabcolsep{0pt}%
    \put(0,0){\includegraphics[width=\unitlength]{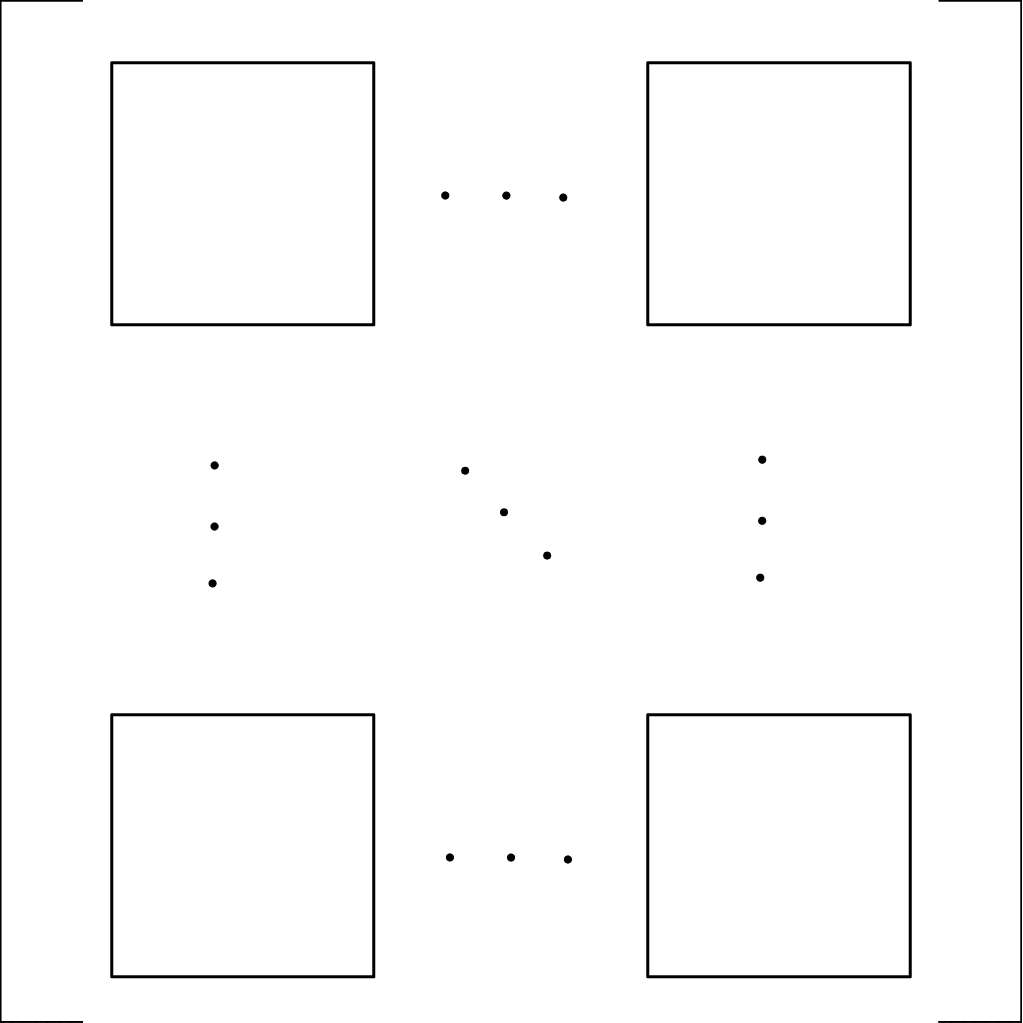}}%
    \put(0.15304585,0.79965451){\color[rgb]{0,0,0}\makebox(0,0)[lt]{\lineheight{1.25}\smash{\begin{tabular}[t]{l}$\mbn{k}^G_{1 1}$\end{tabular}}}}%
    \put(0.66920725,0.79965451){\color[rgb]{0,0,0}\makebox(0,0)[lt]{\lineheight{1.25}\smash{\begin{tabular}[t]{l}$\mbn{k}^G_{1n_{en}}$\end{tabular}}}}%
    \put(0.67034504,0.16071498){\color[rgb]{0,0,0}\makebox(0,0)[lt]{\lineheight{1.25}\smash{\begin{tabular}[t]{l}$\mbn{k}^G_{n_{en}n_{en}}$\end{tabular}}}}%
    \put(0.13838533,0.16168974){\color[rgb]{0,0,0}\makebox(0,0)[lt]{\lineheight{1.25}\smash{\begin{tabular}[t]{l}$\mbn{k}^G_{n_{en}1}$\end{tabular}}}}%
  \end{picture}%
\endgroup%

                \caption{$\mbn{k}^G_{IK}$ to $\mbn{k}_G$}
       \end{subfigure}
	\caption{\centering  Assembly of important parts of element matrices}
	\label{fig_matrix_pattern}
\end{figure}
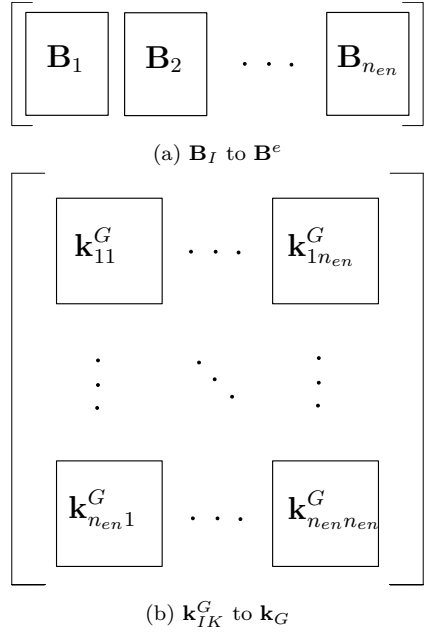

Regarding Eqs. \eqref{lagr_2d}, \eqref{eq:shapefundertransf} and recalling that in GLL quadrature, element nodes coincide with integration points, only one particular shape function $N_I$, and a clearly defined set of derivatives, $N_{J,\alpha}$, $J\in \eta_Q$, are non-zero at the integration point $Q$, which coincides with the node $I$. The shape functions contained in the set $\eta_Q$ all correspond to nodes that lie on the same row or column of point $Q$, as depicted in Fig. \ref{fig_cross_pattern}. In contrast to that, in IGA -- and also in standard $p$-FEM -- all the shape functions of an inner element and their derivatives have non-zero value at each integration point. Therefore, the number of these non-zero value shape functions is equal to $(p+1)^2$ in IGA/$p$-FEM and $(2p+1)$ in SEM methods. The structure of $\mbn{B}^e$ and $\mbn{k}_G$, which contains sub-matrices $\mbn{B}_I$ and $\mbn{k}^G_{IK}$, is shown in Fig. \ref{fig_matrix_pattern}. Considering the number of non-zero shape functions in each method, the number of non-zero sub-matrices in $\mbn{B}^e$ and $\mbn{k}_G$ is equal to $(p+1)^2$ and $(p+1)^4$ for IGA/$p$-FEM and $(2p+1)$ and $(2p+1)^2$ for SEM methods, respectively. By calculating the number of multiplications performed in deriving $\mbn{B}_L$ and $\mbn{k}^G_{IK}$, the total numbers of multiplications in each method, after the numerical integration, are written in Tables~\ref{tab:number_mult_cross} and~\ref{tab:number_mult_standard}.
Please note that computational costs for determining the basis functions and their derivatives are not included, since for IGA this is a recursive process. The influence of these costs is rather negligible in comparison to the costs in the aforementioned tables and can further be reduced by storing and reusing the values within nonlinear iterations. For IGA/$p$-FEM we use the IGA-RMD approach, which requires less effort then the IGA-RMC approach, see~\cite{Dornisch.2015b}.
\begin{table}[h!] 
    \centering 
    \caption{\centering Number of required multiplications for SEM using the cross scheme.} 
    \begin{tabular}{|c|c|c|c|} 
        \hline 
        &Eq. & Total number of multiplications & Asymptotic\\ \hline
        $\mbn{k}_E$ & \eqref{basic-matrices}  & $200(2p+1)^2(p+1)^2+374(2p+1)(p+1)^2$ & $800p^4+\mathcal{O}(p^3)$\\ \hline
        $\mbn{k}_G$ & \eqref{eq:k_g_ik}  & $29(2p+1)^2(p+1)^2+99(2p+1)(p+1)^2$ & $116p^4+\mathcal{O}(p^3)$\\ \hline
    \end{tabular}
    
    \label{tab:number_mult_cross} 
\end{table}

\begin{table}[h!] 
    \centering 
    \caption{\centering Number of required multiplications for IGA-RMD and standard $p$-FEM, or SEM without exploiting the cross scheme.}
    \begin{tabular}{|c|c|c|c|} 
        \hline 
        &Eq. & Total number of multiplications & Asymptotic\\ \hline
        $\mbn{k}_E$ & \eqref{basic-matrices}  & $200(p+1)^6+374(p+1)^4$ & $200p^6+\mathcal{O}(p^5)$\\ \hline
        $\mbn{k}_G$ & \eqref{eq:k_g_ik}  & $29(p+1)^6+99(p+1)^4$ & $29p^6+\mathcal{O}(p^5)$\\ \hline
    \end{tabular}
    \label{tab:number_mult_standard} 
\end{table}
\begin{table}[h!] 
    \centering 
    \caption{\centering Required computations of nodal values for both schemes. Computations can be done before the integration loop as only discrete nodal values occur.} 
    \begin{tabular}{|c|c|c|c|} 
        \hline 
        &Eq. & Total number of multiplications & Asymptotic\\ \hline
        $\mbn{R}_I$ & \eqref{rodrig}  & 
        $47(p+1)^2$ & $47p^2+\mathcal{O}(p^1)$\\ \hline
        $\mbn{T}_{3I}$ & \eqref{eq_T3I}  & $18(p+1)^2$ & $18p^2+\mathcal{O}(p^1)$\\ \hline
        $\mbn{H}_I$ & \eqref{eq_H}  & $21(p+1)^2$ & $21p^2+\mathcal{O}(p^1)$\\ \hline
        $\mbn{T}_I$ & \eqref{var_dir_beta}  & $36(p+1)^2$ & $36 p^2+\mathcal{O}(p^1)$\\ \hline
        $c_3$ & \eqref{eq_MI}  & $4(p+1)^2$ & $4p^2+\mathcal{O}(p^1)$\\ \hline
        $\bar{c}_{10}$ & \eqref{eq_MI}  & $3(p+1)^2$ & $3p^2+\mathcal{O}(p^1)$\\
         \hline
         $c_{11}$ & \eqref{eq_MI}  & $6(p+1)^2$ & $6p^2+\mathcal{O}(p^1)$\\
         \hline\hline
         $\sum$ & & $135(p+1)^2$ & $135p^2+\mathcal{O}(p^1)$\\ \hline
    \end{tabular}
    
     \label{tab:number_mult_nodal}
\end{table}
\begin{remark}
    the multiplications in $\mbn{k}_E$ are multiplications of large matrices, which are optimized very well in modern compilers. In $\mbn{k}_G$, there are many small multiplications which have to be called one by one. We did not include symmetry in our considerations. For $\mbn{k}_E$, this reduces the number of multiplications asymptotically by one half. The lower right part of  $\mbn{k}_G$, see Eq.~\eqref{eq:k_g_ik}, is already symmetric and constitutes the most expensive part of $\mbn{k}_G$. Thus, the computation of $\mbn{k}_G$ is not reduced significantly by exploiting symmetry.
\end{remark}
\par In deriving Tables \ref{tab:number_mult_cross}, \ref{tab:number_mult_standard} and \ref{tab:number_mult_nodal}, it is assumed that division is approximately equal to multiplication. Additions are neglected since they typically require significantly less numerical effort. Additionally, the $\sin$ and $\cos$ functions of the nodal vector of rotation $\omega_I$ are precomputed and stored outside the integration loops. The stored values are simply called during the calculation of $\mbn{R}_I$, $\mbn{H}_I$ and $\mbn{M}_I$. Thus, $\sin\omega_I$ and $\cos\omega_I$ each has to be called only $(p+1)^2$ times per element. The effort for computing one $\sin$ or $\cos$ function is in the range of a few multiplications. Thus, the total effort for the angular functions is in the range of $10(p+1)^2$ to $20(p+1)^2$, which aligns with the values of Table~\ref{tab:number_mult_nodal}. In counting the multiplications in Tables \ref{tab:number_mult_cross}, \ref{tab:number_mult_standard} and \ref{tab:number_mult_nodal}, we assumed that the results of all multiplications which occur several times are stored and reused.

\begin{figure}[ht] 
	\center{\includegraphics[scale=0.9]{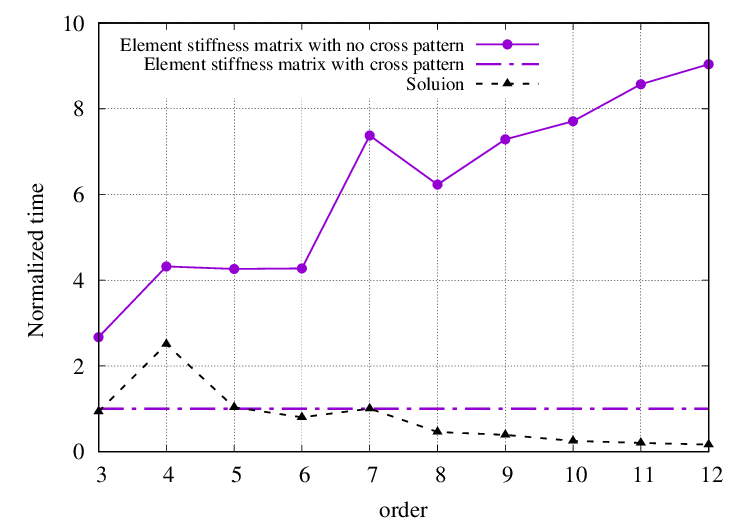} }
	\caption{\centering Normalized time of calculation of element stiffness matrix without using cross pattern and normalized time solution vs order of element }\label{fig_time_cross}
\end{figure}

From Tables~\ref{tab:number_mult_cross} and~\ref{tab:number_mult_standard}, it is evident that the computational cost to calculate the stiffness matrix of an element is significantly reduced for SEM methods compared to IGA/$p$-FEM methods, especially for higher-order elements. To better illustrate this fact, efficiency of utilizing the cross pattern is illustrated in Fig.~\ref{fig_time_cross}. In this figure, a single SEMI element is used to model the Scordelis-Lo roof. While the material properties and dimensions of this example will be detailed in the next section, our current focus is on evaluating the efficiency of using the cross pattern. To generate Fig.~\ref{fig_time_cross}, a single load step is performed, and the computation times for both calculation of element stiffness matrix without exploiting the cross scheme and solution of the linear system of equations are normalized against the stiffness matrix calculation time, when using the cross pattern. In theory, the solution time should remain the same whether the cross pattern is applied or not, as the resulting stiffness matrices are identical in both cases. An asymptotic trend emerges in very high-order elements, where the polynomial order $p$ dominates the numerical coefficients in Tables~\ref{tab:number_mult_cross} and \ref{tab:number_mult_standard}. Thus, the efficiency gains from neglecting zero entries in $\mbn{B}_L$ and $\mbn{k}^G_{IK}$ are evident, highlighting that for higher-order elements, the predominant computational cost lies in constructing the stiffness matrix. We anticipate that for an extremely fine mesh, the cost of solving the system may increase non-linearly, yet the computational burden of generating element stiffness matrices will also grow dramatically. Thus, Fig.~\ref{fig_time_cross} serves as a valuable reference for understanding the significance of the cross pattern in computational efficiency. Finally, it should be noted that such a significant increase in speed is achieved simply by adding a few lines of code that check whether the integration point lies on the cross pattern or not.

Another issue that arises when higher-order elements are used is the high number of internal nodes, which have no connectivity with other elements. These nodes increase the computational cost of solving the equations and enlarge the global stiffness matrix. In \cite{payette_2014}, static condensation is employed to exclude these DOFs, and it is concluded that the method is highly effective and efficient. However, we do not include static condensation in our considerations and postpone this to further works.

What has been presented so far in this section is the comparison of deriving  the stiffness matrices of only one high order element in SEM and IGA and it is shown that SEM needs much less computations even in comparison with the simplest IGA formulation, namely IGA-RMD. On the other hand, the IGA methods benefit from $k$-refinement, which allows the use of higher-order elements with a finer mesh while keeping the number of DOFs relatively low. However, since industrial CAD files often include a significant number of internal knots even in the coarsest mesh, the optimal results of $k$-refinement cannot be fully achieved. Thus, computations using SEM can remain competitive in practical applications, even when $k$-refinement is employed in IGA.

\section{Numerical Examples} \label{sec_4}
This section presents four geometrically nonlinear numerical examples which are solved using the newly presented SEMI approach. To assess accuracy and efficiency, we provide also results using the SEMN approach, initially proposed for linear analysis in~\cite{nima_2024}, and using an IGA shell formulation, proposed in~\cite{Dornisch.2016}. Two of these examples are nonlinear versions of the well-known shell obstacle course proposed in~\cite{belyt_1985}. They are chosen due to their susceptibility to shear or extensional locking. We maintain the same structural dimensions and material properties as in \cite{belyt_1985}, but for one of them we increase the applied loads to induce large deformations and rotations of the structure. The two further numerical examples are free form surfaces, which aim at assessing the ability of SEMI to capture geometries with large curvature gradients. The first one is taken from~\cite{Dornisch.2013}, while the second one is newly proposed by us and contains a very steep change of curvature to showcase the behavior of all formulations in extreme cases. Computations with all other methods shown here - namely SEMN, IGA-RMD and IGA-RMC - make use of the exact geometry description by NURBS, while SEMI approximates the geometry. Thus, these examples are able to further elaborate on the question if the exact representation of the geometry is necessary to provide high accuracy results. By comparing 
\begin{itemize}
    \item SEMN: exact geometry by using NURBS geometry description for position vectors, director vectors and Jacobian. Derivatives of unknown fields approximated by Lagrange shape functions
    \item SEMI: exact position vectors in the integration points (nodal values computed from the NURBS geometry); director vector, all derivatives and Jacobian approximated by Lagrange shape functions
\end{itemize}
we can clearly assess the influence of the exact geometry description, which has been shown in~\cite{Oesterle.2022} to be of minor influence, but is revisited here within a different setting. The further comparison to
\begin{itemize}
    \item IGA-RMD: all quantities approximated by NURBS basis functions, discrete rotational formulation which tends to artificial thinning
    \item IGA-RMC: all quantities approximated by NURBS basis functions, continuous rotational formulation which ensures correct shell thickness
\end{itemize}
allows assessing multiple points such as influence of high-continuity, influence of rotational formulation and influence of exact representation of geometry and derivatives.

However, the main intention of this paper is to show the efficiency and accuracy of the newly proposed SEMI shell. We mainly want to use $p$-refinement with orders up to $p=15$ to get rid of locking effects while maintaining reasonable computational costs by using the cross scheme. In this setting, the number of DOFs between SEM and IGA methods is the same, and IGA cannot profit from the higher continuity. Thus, we also added combined $h$-refinement and order elevation to get a fair comparison. In~\cite{nima_2024} we already assessed the condition number of SEM methods, which has been shown to be growing much slower with rising order than for IGA methods. Thus, we can use really high orders without running into numerical problems. While in the linear cases in~\cite{nima_2024}, IGA has shown appropriate results up to $p=13$ despite bad condition numbers, we will show that for the nonlinear examples IGA fails for significantly lower orders while SEMI works reliably. Due to the high numerical effort to compute condition numbers, we do not provide them for our numerical examples but refer to the linear computations in~\cite{nima_2024}.
\subsection{Scordelis-Lo Roof}
\label{sec:scordelis}
\begin{figure}[ht] 
 	\center{\includegraphics[scale=0.7]{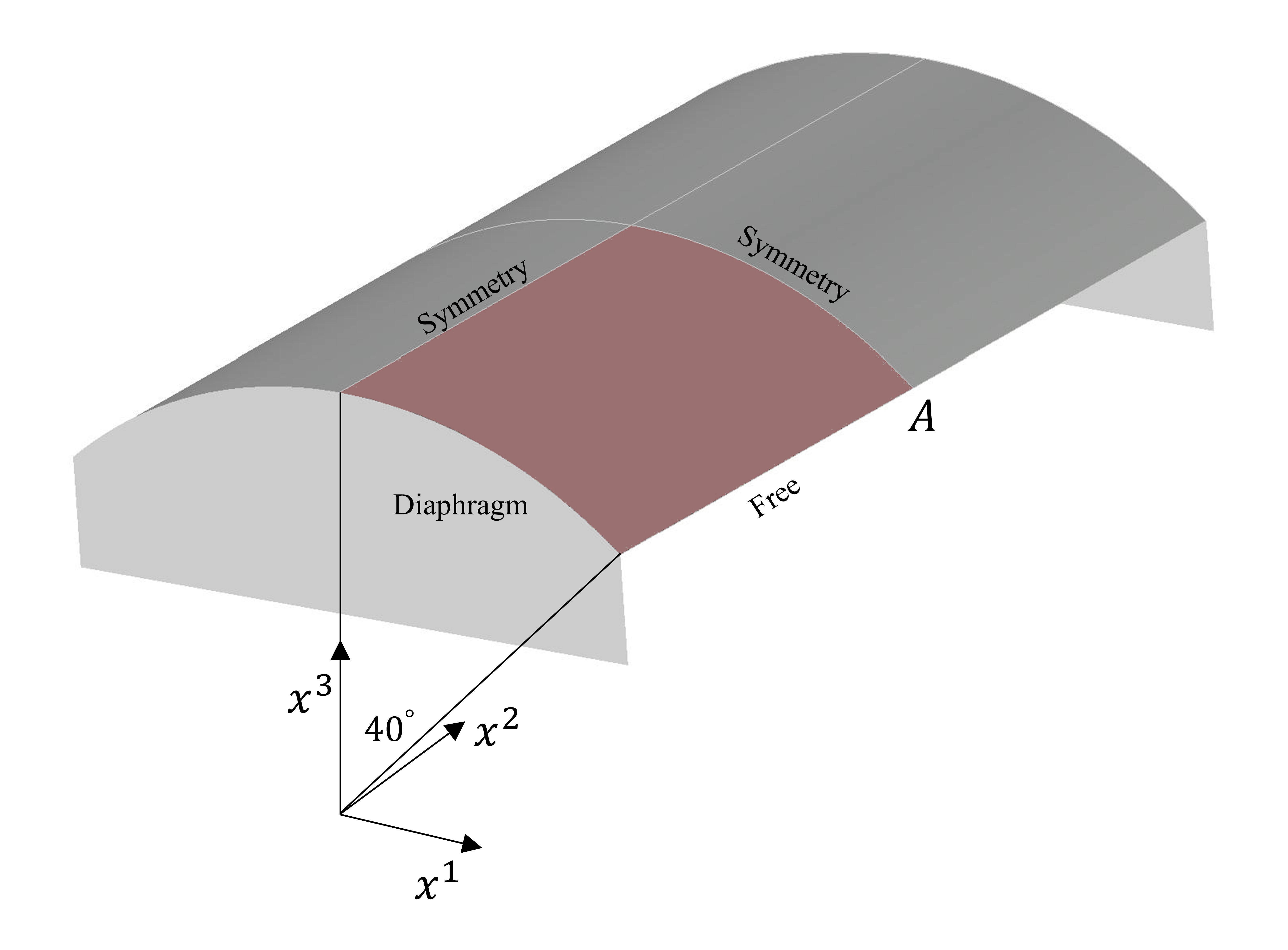} }
\caption{\centering Schematic for Scordelis-Lo roof}\label{fig_scordelis}
\end{figure}

In this example large deformations of the Scordelis–Lo roof are analyzed. At each end of the structure there is a rigid diaphragm, which is parallel to the \mbox{$x^1$-$x^3$} plane and restrains the in-plane displacements in $x^1$ and $x^3$ directions and prevents in-plane rotations. The roof is subjected to a distributed uniform gravity load equal to $90$ per unit area. The radius of the roof is $25$ and the thickness of the shell is $0.25$. The length of the whole roof is $50$. Young’s modulus and Poisson’s ratio are $4.32 \cdot 10^8$ and $0$, respectively. Considering symmetry, only a quarter of the shell, as shown in Fig. \ref{fig_scordelis}, is modeled, and the vertical displacement at point A is compared using different methods.  The load is applied in one step.

\begin{figure}[ht] 
	\center{\includegraphics[scale=0.9]{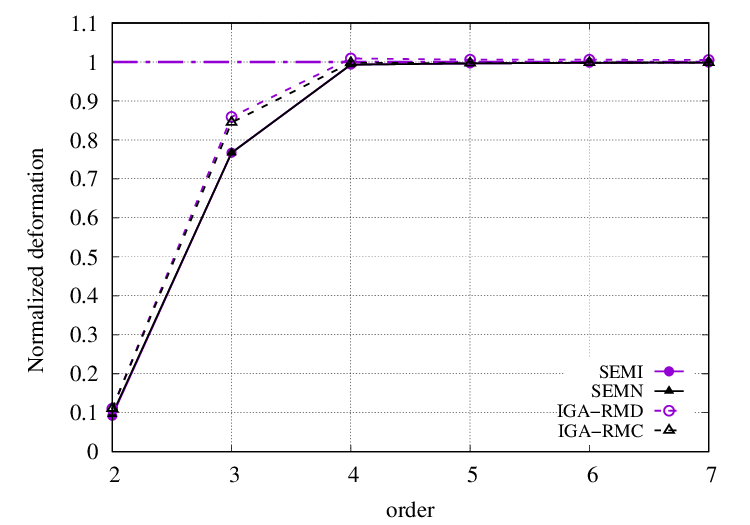} }
	\caption{\centering Scordelis-Lo roof: Normalized deformation for \textit{p}-refinement  vs order of element}
	\label{fig_sc_pref_1}
\end{figure}

\begin{figure}[ht] 
	\center{\includegraphics[scale=0.9]{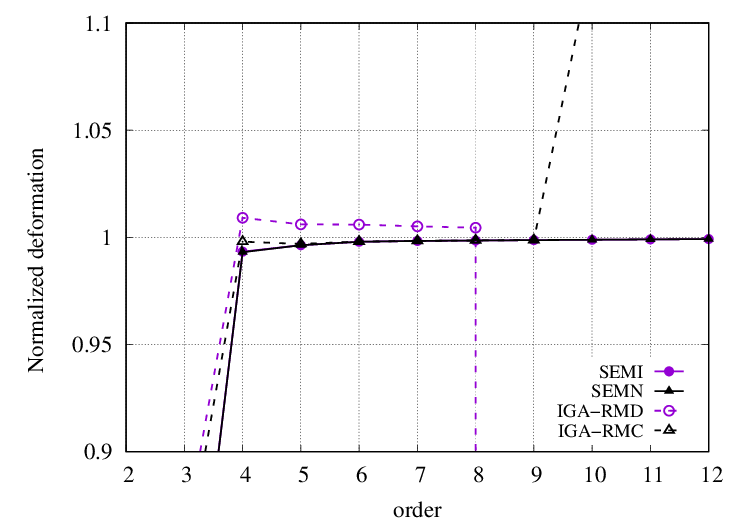} }
	\caption{\centering Scordelis-Lo roof: Normalized deformation for higher order \textit{p}-refined elements  vs order of element}
	\label{fig_sc_pref_2}
\end{figure}

\begin{figure}[ht]
	\begin{subfigure}[t]{\linewidth}
			\includegraphics[width=120mm]{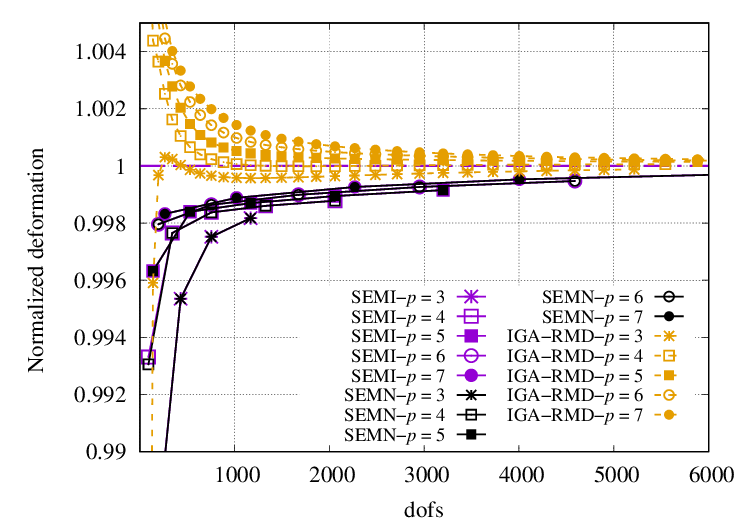} 
                \caption{comparison to IGA-RMD}
                \label{fig_sc_3_RMD}
       \end{subfigure}
       \begin{subfigure}[t]{\linewidth}
                \includegraphics[width=120mm]{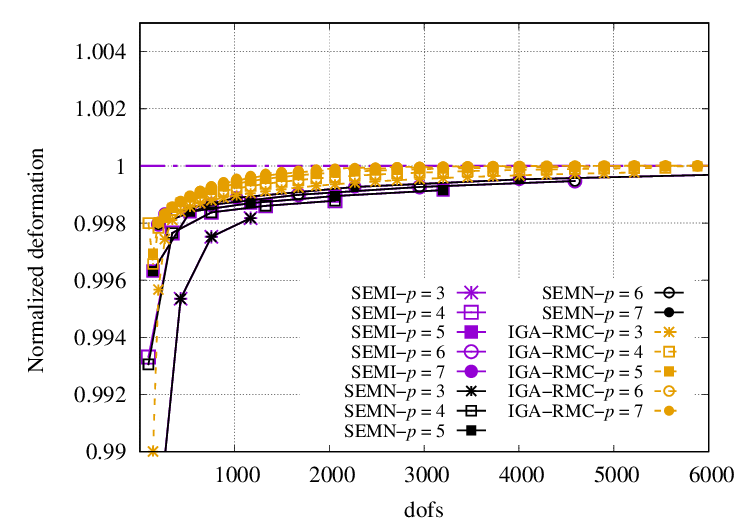}
                \caption{comparison to IGA-RMC}
                \label{fig_sc_3_RMC}
       \end{subfigure}
	\caption{\centering Scordelis-Lo roof: Normalized deformation vs DOFs for \textit{h}-refinement of SEMI and SEMN, compared to IGA}
	\label{fig_sc_3}
\end{figure}

\begin{figure}[ht]
	\begin{subfigure}[t]{\linewidth}
			\includegraphics[width=120mm]{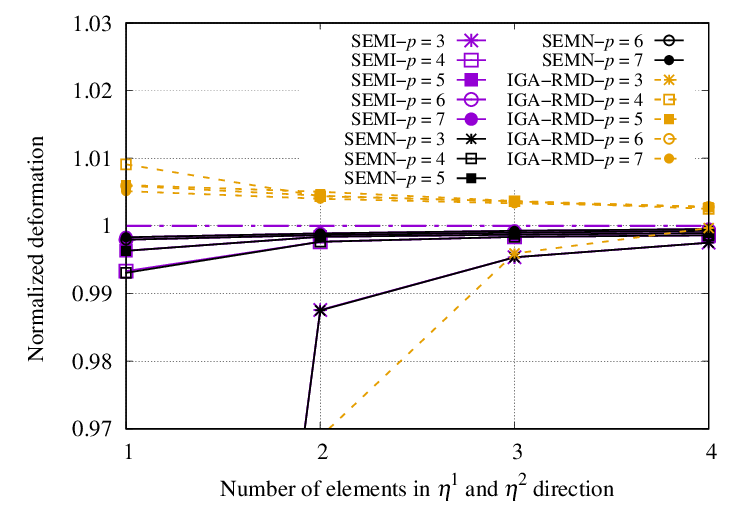} 
                \caption{comparison to IGA-RMD}
                \label{fig_sc_4_RMD}
       \end{subfigure}
       \begin{subfigure}[t]{\linewidth}
                \includegraphics[width=120mm]{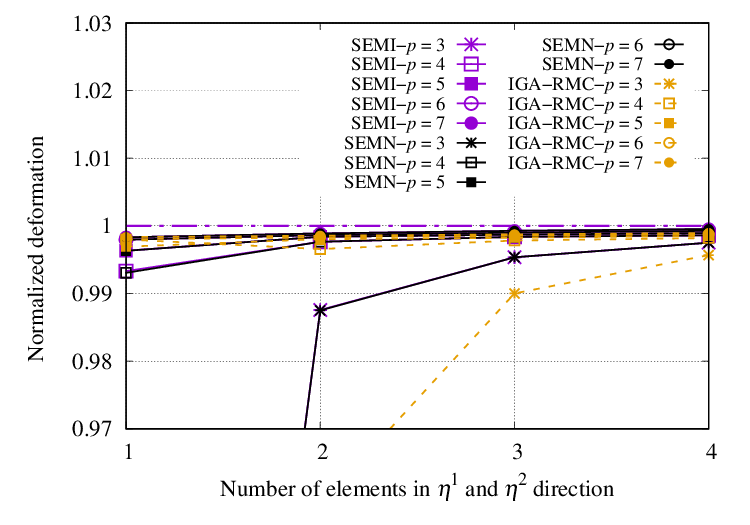}
                \caption{comparison to IGA-RMC}
                \label{fig_sc_4_RMC}
       \end{subfigure}
	\caption{\centering Scordelis-Lo roof: Normalized deformation vs number of elements  for \textit{h}-refinement of SEMI and SEMN, compared to IGA}
	\label{fig_sc_4}
\end{figure}

Fig.~\ref{fig_sc_pref_1} shows a comparison between different methods using $p$-refinement, i.e. one element with rising orders of the shape functions. The deformation is normalized by the converged result $u^3_A=-0.25356483$ of IGA-RMC with a very fine mesh ($40\times 40$ elements with $p=8$). As observed in Fig.~\ref{fig_sc_pref_1}, at lower orders, IGA methods provide slightly better results, but beginning with $p=4$, all methods yield almost the same result and are accurate enough for practical applications. It is to be noted, that SEMN and SEMI converge from below to the same result as IGA-RMC, which is the more accurate IGA formulation. IGA-RMD, on the other hand, yields slightly higher deformations and does not converge properly for pure $p$-refinement for higher-order elements. To investigate the behavior of the elements at very high orders, the normalized deformation for orders up to $p=12$ is presented in Fig.~\ref{fig_sc_pref_2} using a smaller scale around the exact solution. IGA-RMD starts to diverge for $p=9$, IGA-RMC for $p=10$, while SEMI and SEMN remain numerically stable and converge to the correct deformation. The unstable behavior of the IGA elements is evident and, as discussed regarding the increase in the condition number in the IGA method, also predictable. The fact that IGA-RMC yields slightly better numerical stability than IGA-RMD can also be seen in~\cite{nima_2024} and is attributed to its more accurate rotational formulation.
The normalized deformation for $h$-refinement is plotted in Fig.~\ref{fig_sc_3} and  Fig.~\ref{fig_sc_4} versus the number of DOFs and the number of elements, respectively. Order $p=2$ is excluded from both figures because the results fall outside the depicted ranges due to severe locking. The figures clearly show that there is no significant difference between SEMN and SEMI. Thus, the exact representation of the geometry of SEMN does not yield any positive effect. The accuracy of the IGA-RMD computations in Fig.~\ref{fig_sc_3_RMD} is higher than for SEMI for low orders, but deteriorates with rising order, which is line with the findings in~\cite{Dornisch.2013}. This clearly shows that even though both SEMI and IGA-RMD use the same rotational formulation, IGA-RMD is not suitable for high orders while SEMI converges correctly for rising orders. The accuracy of IGA-RMC computations in Fig.~\ref{fig_sc_3_RMC} is always higher than for SEMI, and the convergence behavior is correct. Two points have to be noted for this comparison. Firstly, the rotation formulation of IGA-RMC is significantly more complex than of SEMI. Secondly, IGA yields significantly less DOFs per element due to the higher continuity. Thus, we also compare SEMI and SEMN to IGA with number of elements used as abscissa. From Fig.~\ref{fig_sc_4_RMD}, it can be deduced that SEMI and SEMN produce significantly better results than IGA-RMD with the same mesh. The situation is different for IGA-RMC, which is depicted in Fig.~\ref{fig_sc_4_RMC}. If using one element only, IGA-RMC is more accurate than SEMI for low orders and as accurate as SEMI for high orders. But beginning from two elements per edge, SEMI has a slightly higher accuracy per number of element than IGA-RMC for all orders. This shows again that the approximation of the geometry of SEMI is no real drawback.

\subsection{Hemisphere with a hole}
\begin{figure}[ht] 
	\center{\includegraphics[scale=0.7]{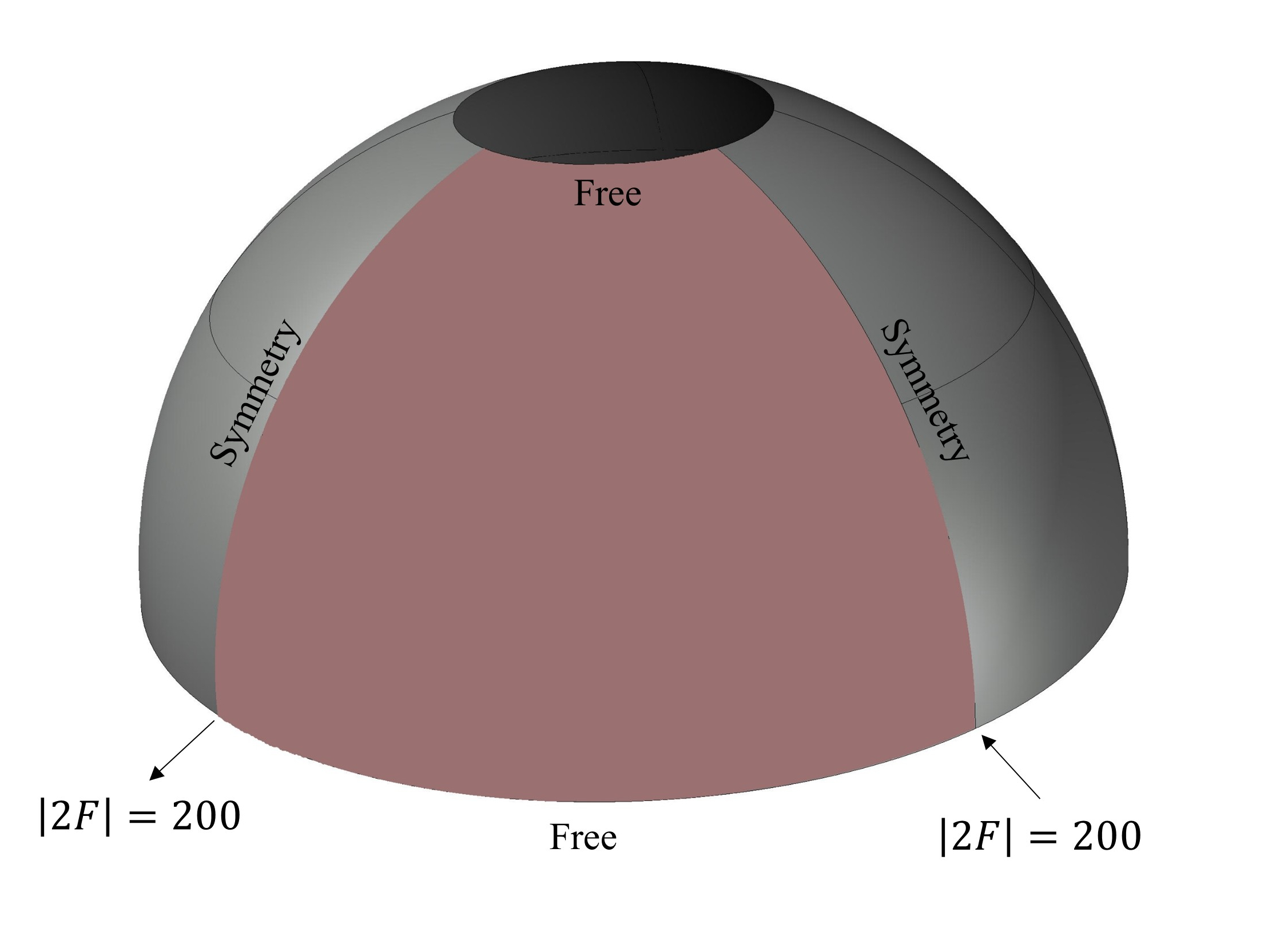} }
	\caption{\centering Schematic for hemispherical shell with hole}\label{fig_hemisphere_0}
\end{figure}

In Fig.~\ref{fig_hemisphere_0}, a hemispherical shell with a $18^\circ$ hole on top-center and a radius of $10$ is depicted.Young's modulus and Poisson's ratio are $6.825\cdot 10^7$ and $0.3$, respectively. The thickness of the shell is $0.04$. Radial forces $2F$ alternating at $90^\circ$ are applied to the shell as shown. Regarding the symmetry, only a quarter of the hemisphere subjected to load $F$ is analyzed. The radial deformation under the point load is investigated. In the original problem the value $\left| 2 F \right| = 2$ is chosen \cite{belyt_1985}. However, to yield significantly larger rotations and deformations in this study, a value of $\left| 2 F \right| = 200$ is used. The deformation is normalized by the radial deformation $u=5.86799$, which has been computed with $200\times 200$ elements of order $p=6$ using the IGA-RMC formulation.
\par
\begin{figure}[ht] 
	\center{\includegraphics[scale=0.9]{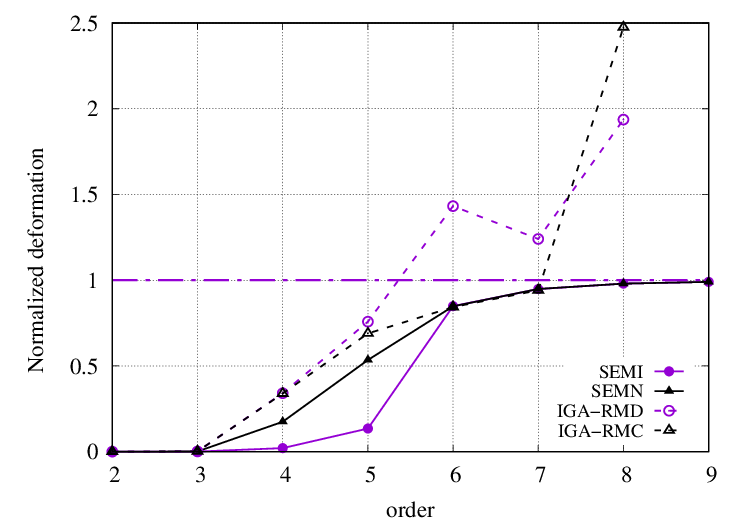} }
	\caption{\centering Hemispherical shell: Normalized deformation for \textit{p}-refinement  vs order of element}
	\label{fig_hem_pref}
\end{figure}

\begin{figure}[ht]
	\begin{subfigure}[t]{\linewidth}
			\includegraphics[width=120mm]{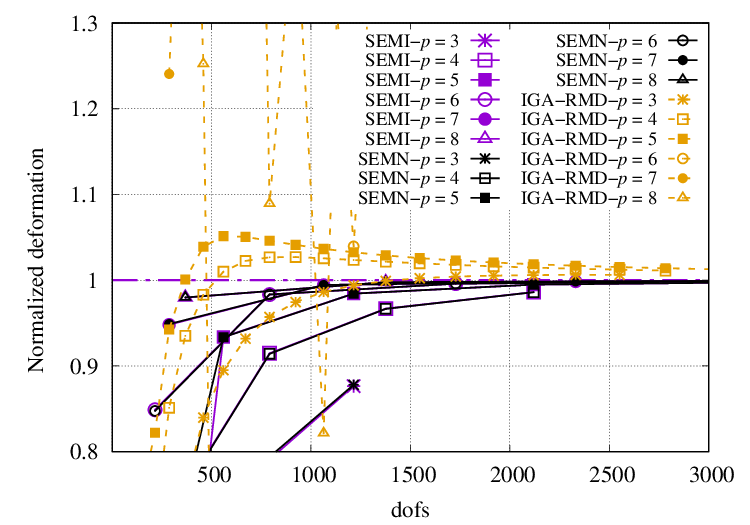} 
			\caption{comparison to IGA-RMD}
        \end{subfigure}
        \begin{subfigure}[t]{\linewidth}
			\includegraphics[width=120mm]{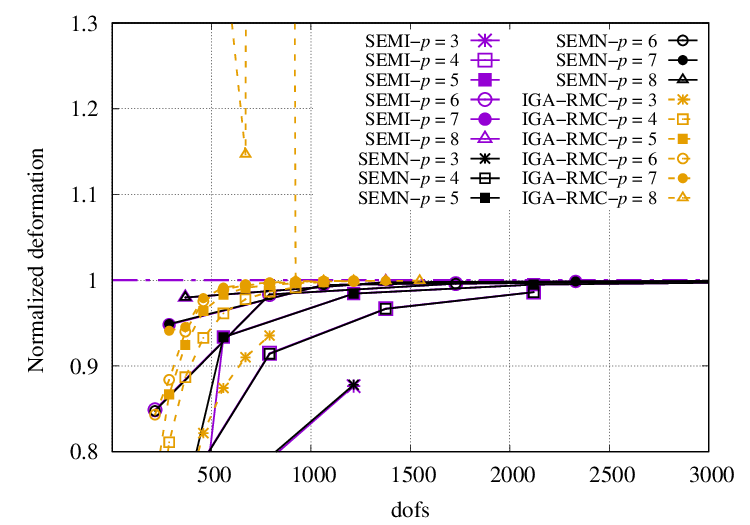}
			\caption{comparison to IGA-RMC}
	\end{subfigure}
	\caption{\centering Hemispherical shell: Normalized deformation vs DOFs for \textit{h}-refinement of SEMI and SEMN, compared to IGA}
	\label{fig_hem_href}
\end{figure}

 Fig.~\ref{fig_hem_pref} demonstrates the superiority of SEM-type methods over IGA in $p$-refinement while only one element is used. For orders $p=2$ to $p=5$, the IGA formulations and SEMN profit from the exact geometry and significantly outperform SEMI. However, the error level of all formulations is still far from acceptable. Beginning from $p=6$, the exact representation of the geometry seems to have no influence anymore since SEMI, SEMN and IGA-RMC coincide very well. This clearly shows that the solution approximation error is way more important than the geometry approximation error, which is in line with the findings of~\cite{Oesterle.2022}. For rising orders, the SEMI method converges correctly and is stable, while IGA-RMD and IGA-RMC start to get numerically unstable and diverge completely for $p=6$ and $p=8$, respectively. In the Scordelis-Lo roof in Sec.~\ref{sec:scordelis}, a single-curved shell example, IGA-RMD and IGA-RMC start to diverge for orders $p=9$ and $p=10$, respectively. This is a first hint that despite the use of the the exact geometry in IGA, more complex geometries should be computed using lower orders of basis functions. This hypothesis will be studied in the further examples. Recalling that SEMI uses the same rotational formulation as IGA-RMD, a very interesting point can be seen in Figs.~\ref{fig_hem_pref} and~\ref{fig_hem_href}. SEMI converges from below in the same manner as the more complex IGA-RMC formulation, with visually the same accuracy for $p=6$ and $p=7$ for one element. This is due to the Kronecker-delta property of SEMI, which prevents artificial thinning as effective as IGA-RMC. Thus, it can be clearly attested that for SEMI, the simpler discrete rotational formulation proposed in this paper is entirely sufficient.   

\begin{figure}[ht]
	\begin{subfigure}[t]{\linewidth}
			\includegraphics[width=120mm]{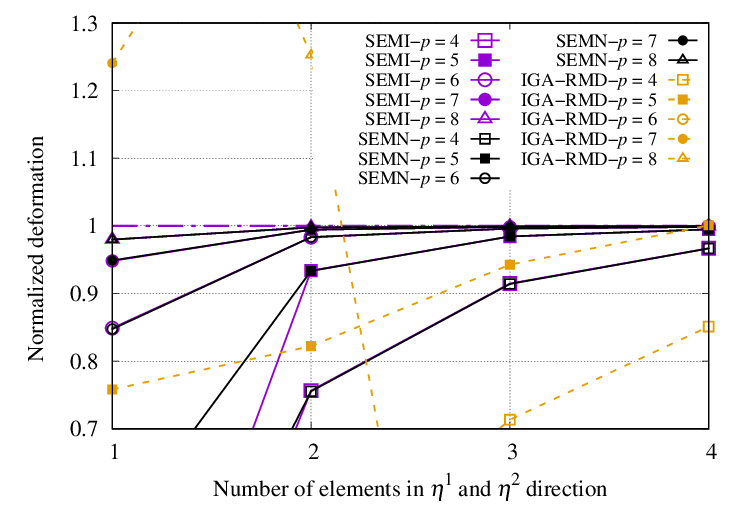} 
			\caption{comparison to IGA-RMD}
        \end{subfigure}
        \begin{subfigure}[t]{\linewidth}
			\includegraphics[width=120mm]{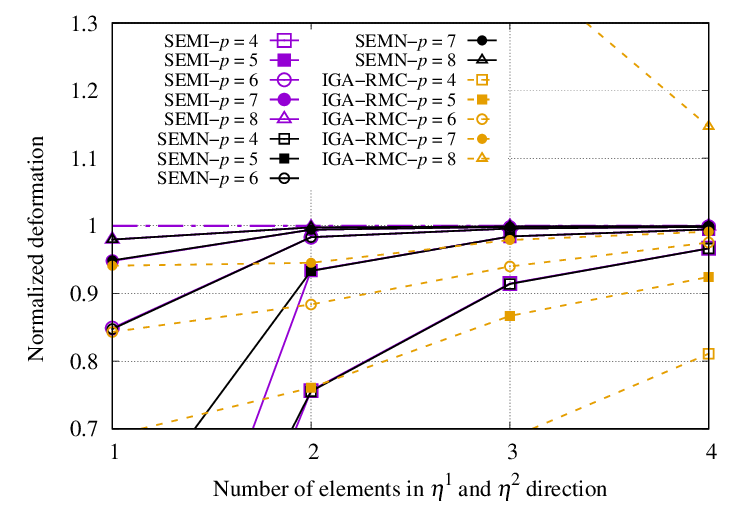}
			\caption{comparison to IGA-RMC}
	\end{subfigure}
	\caption{\centering Hemispherical shell: Normalized deformation vs number of elements  for \textit{h}-refinement of SEMI and SEMN,  compared to IGA}
	\label{fig_hem_href_elem}
\end{figure}

\par In Fig.~\ref{fig_hem_href}, a comparison between the different methods for $h$-refinement is depicted. Results for elements of low order are excluded from Figs.~\ref{fig_hem_href} and \ref{fig_hem_href_elem} if they fall outside the depicted ranges. For the IGA methods, $k$-refinement is applied. Thus, for the same number of elements, the number of DOFs is generally smaller in IGA compared to SEMI or SEMN. As shown in Fig.~\ref{fig_hem_href}a, SEMI and SEMN generally provide better results for orders above 5, while at lower orders, IGA-RMD achieves more accurate results. Thus, for the same number of DOFs, an optimal choice would be to use higher-order elements with fewer elements in SEMI or SEMN, whereas in IGA-RMD, lower-order elements with a greater number of elements yield better accuracy. Figure~\ref{fig_hem_href}b shows a significant improvement in the results of IGA-RMC compared to IGA-RMD, owing to a more accurate rotational formulation. However, it also demonstrates that by selecting high-order elements with a small number of elements in SEMI and SEMN, we can achieve results close to those of IGA-RMC, despite the latter's more sophisticated formulation. Although, as explained, at a specific number of DOFs, IGA-RMC provide better results in comparison with SEMN and SEMI, if we consider the number of elements, SEMI and SEMN have the upper hand as depicted in Fig.~\ref{fig_hem_href_elem}. This point can be considered as an advantage if complicated material laws, such as plasticity, are used, because then the number of elements has a significant influence on the total computational costs. Furthermore, we can see again that the approximation of the geometry of SEMI is no drawback, since already for $2\times 2$ elements of order $p=8$, the SEMI solution cannot be distinguished by the eye from the reference solution. Interestingly, by going from one to four elements, SEMI significantly improves its accuracy, while in IGA-RMC there is only little gain in accuracy for this refinement step.

\subsection{Double curved free form geometry}
\begin{figure}[ht] 
	\center{\includegraphics[scale=0.5]{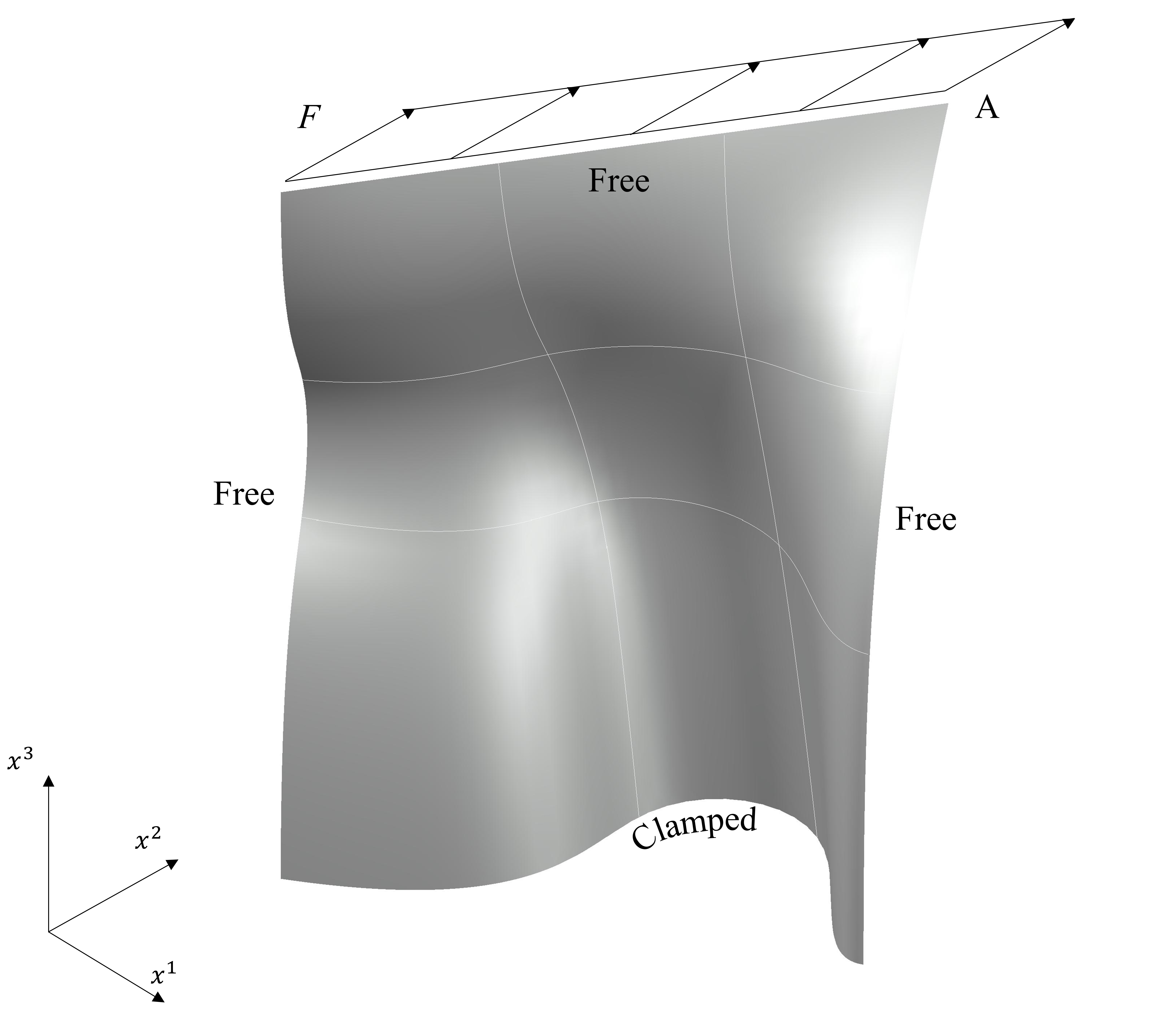} }
	\caption{\centering Schematic for the double curved free form geometry}\label{fig_ff_0}
\end{figure}

\begin{figure}[ht] 
	\center{\includegraphics[scale=0.9]{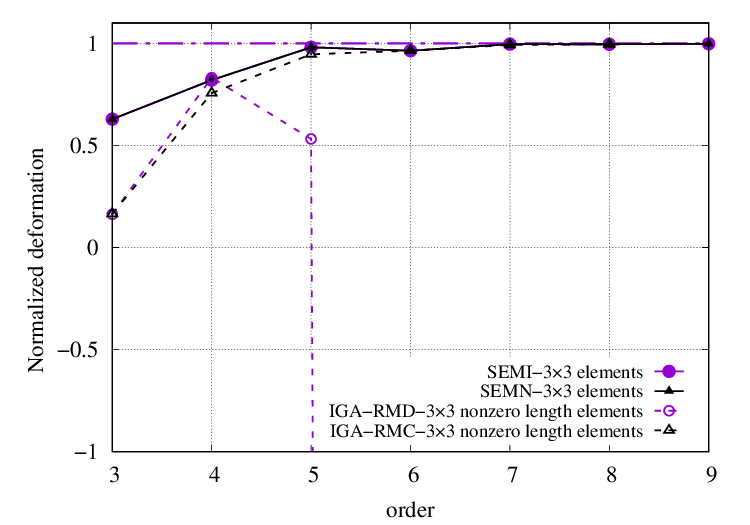} }
	\caption{\centering Free form shell: Normalized deformation for \textit{p}-refinement  vs order of element for the SEM cases with evenly spaced nodes in SEM parametric space}
	\label{fig_ff_pref}
\end{figure}

\begin{figure}[ht]
	\begin{subfigure}[t]{\linewidth}%
		\includegraphics[width=120mm]{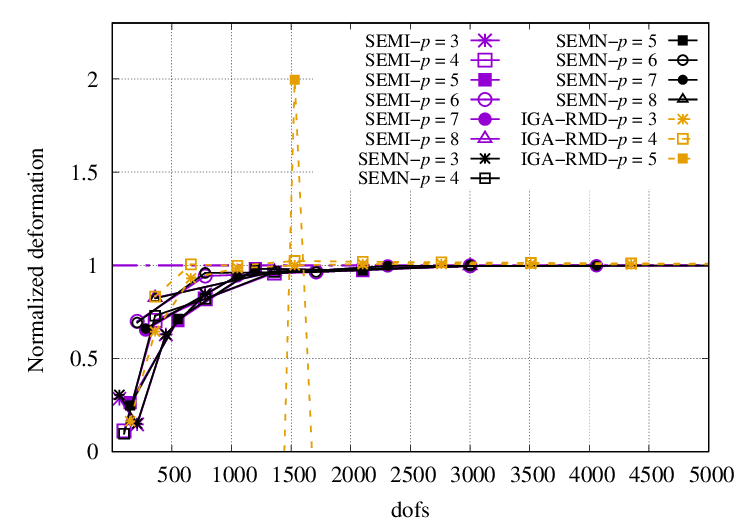}
		\caption{comparison to IGA-RMD}
        \end{subfigure}
        \begin{subfigure}[t]{\linewidth}
				\includegraphics[width=120mm]{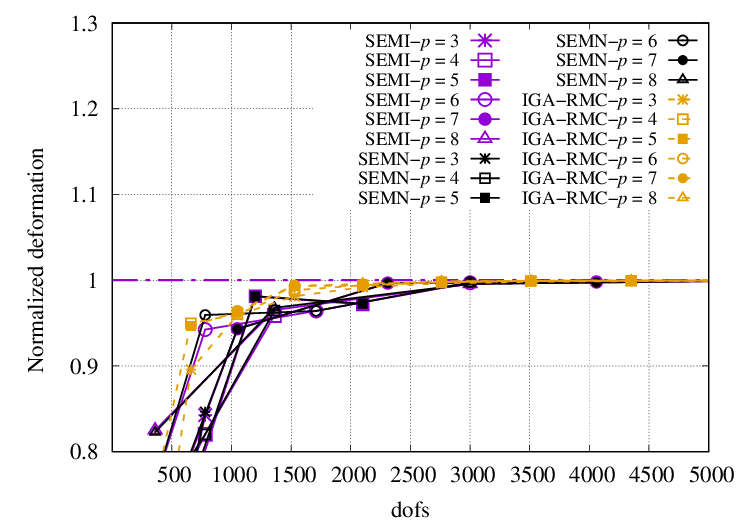}
				\caption{comparison to IGA-RMC}\label{fig_ff_href_b}
	\end{subfigure}
	
	\caption{\centering Free form shell: Normalized deformation vs DOFs for \textit{h}-refinement of SEMI and SEMN with evenly spaced nodes in SEM parametric space, compared to IGA}
	\label{fig_ff_href}
\end{figure}

\begin{figure}[ht] 
	\center{\includegraphics[scale=0.9]{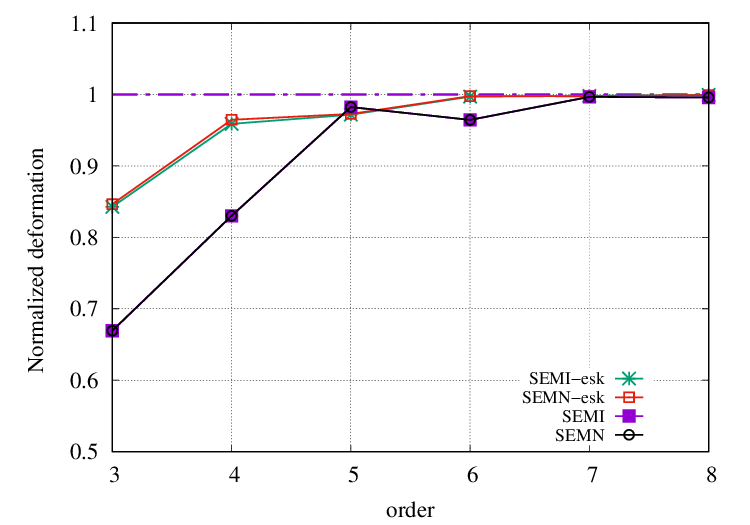} }
	\caption{\centering Free form shell: Normalized deformation for $4\times4$ mesh vs order of element for the SEM scenarios, 'esk' refers to the evenly spaced knots, i.e. the first scenario}
	\label{fig_ff_4x4}
\end{figure}

This example focuses on the comparison of the different explained methods to deal with complex surfaces and comparing the results of different refinement methods. Therefore, a double curved free from surface depicted in Fig.~\ref{fig_ff_0}  with the thickness of 0.1 is chosen to analyze. The geometry of the surface can be constructed by using the knot vector and control points provided in  \ref{app2}; see also~\cite{Dornisch.2013}. A uniform line load of $F=10$ per unit length in $x^2$ direction is applied at the top edge of the surface. Young's modulus and Poisson's ratio are $1.2\cdot10^6$ and 0.3, respectively. The converged displacement at point A in direction of the $x^2$ axis utilized to normalize the result is 0.734541. There are two internal knots in the knot vectors of both parametric directions in the CAD model. Hence, the coarsest mesh in IGA consists of $3\times 3$ elements. In the $p$-refinement study for SEM, we have used a mesh which coincides with those non-zero length elements. The elementation for the IGA computations is defined in the CAD file. Additionally, order elevation is applied. The results are depicted in Fig. \ref{fig_ff_pref}. SEMI provides better results than IGA-RMC in lower orders despite the simpler rotation interpolation and the inexact geometry. The exact geometry of SEMN brings no advantage over SEMI.  IGA-RMD diverges already at order $p=5$. 
The comparison for $h$-refinement is provided in Fig.~\ref{fig_ff_href}. IGA-RMD yields only stable results up to an order $p=4$. The accuracy of IGA-RMC is slightly higher than of SEMI, mainly due to the higher number of elements. Thus, IGA benefits from $k$-refinement. In IGA all three initial knot spans are subdivided equally for the refinement study, yielding $9, 36, 81, 144,\ldots$ elements. For SEM, the most refined mesh we computed for this example consists of $4\times 4$ elements. We have defined two scenarios to be utilized in calculation of the results. In the first case, SEM elements are evenly spaced in the IGA parametric domain (abbreviated as 'esk' in Fig.~\ref{fig_ff_4x4}). In the second scenario, the IGA knots are used as initial element boundaries. For further refinement, elements are subdivided. For $4\times 4$ elements, this results in nodal values $\eta^{\alpha}_i=0,\frac{1}{3},\frac{1}{2},\frac{2}{3},1 $. 
			
				
				
				

\begin{figure}[ht] 
	\center{\includegraphics[scale= 0.9]{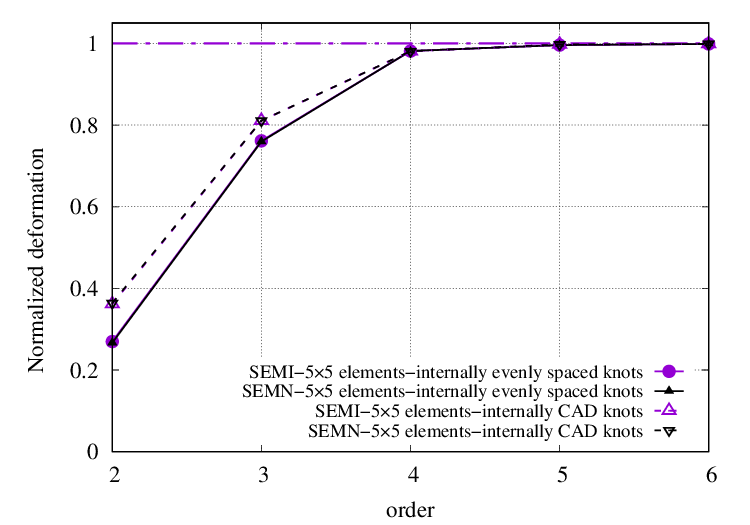} }
	\caption{\centering Free form shell: Normalized deformation for \textit{p}-refinement  vs order of element for SEMI and SEMN, 5$\times$5 mesh with evenly spaced and as per CAD internal nodes}\label{fig_ff_cad_even}
\end{figure}

The aim of these scenarios is to investigate the effect of using different meshing methods, which are described in Section \ref{section_SEMN}. From intuition, we assume that if we use the CAD internal knots as the border of elements for meshing, regarding the reduction of continuity at internal knots, it will be a more 
 logical choice and provide more accurate results in SEMI and SEMN. Fig.~\ref{fig_ff_href}, which is based on the first scenario, shows that the results of IGA methods are a little more accurate in coarser meshes but in finer meshes, SEM methods yield a similar accuracy as IGA-RMC. Fig. \ref{fig_ff_4x4} shows the comparison between the two meshing scenarios. Contrary to our expectations, evenly spaced elements ('esk') provides more accurate results in lower orders. This may be due to the fact that we have depicted the normalized displacement of point A which is directly affected by the size of the side elements. If $4\times 4$ elements of equal length are used, the side elements are smaller than the side elements in the second scenario. To further investigate this, the results for a $5\times5$ mesh is shown in Fig. \ref{fig_ff_cad_even}. 
 For the CAD-oriented mesh we use the CAD knots and refine the side elements. So, coordinates of the border of elements in IGA space are 
\begin{equation}\label{ff_5x5_CAD}
	\mbn{\Xi}=\left[0,\frac{1}{6},\frac{1}{3},\frac{2}{3},\frac{5}{6},1 \right]\,.
\end{equation}
 To allow judging about the necessity of setting element boundaries at CAD knots, we now put two element boundaries at $\frac{1}{6}$ and $\frac{5}{6}$. The additional two inner knots are obtained by equal subdivision between the aforementioned two knots. Accordingly, the coordinates of the borders for the internally evenly spaced meshing in Fig.~\ref{fig_ff_cad_even} is  
\begin{equation}\label{ff_5x5_evenly}
	\mbn{\Xi}=\left[0,\frac{1}{6},\frac{7}{18},\frac{11}{18},\frac{5}{6},1 \right]\,.
\end{equation} 
 
Fig. \ref{fig_ff_cad_even} clearly shows that for lower order elements, if we use CAD borders, we gain slightly better results.

\subsection{Double curved free form NURBS based geometry}
\begin{figure}[ht] 
	\center{\includegraphics[scale=0.75]{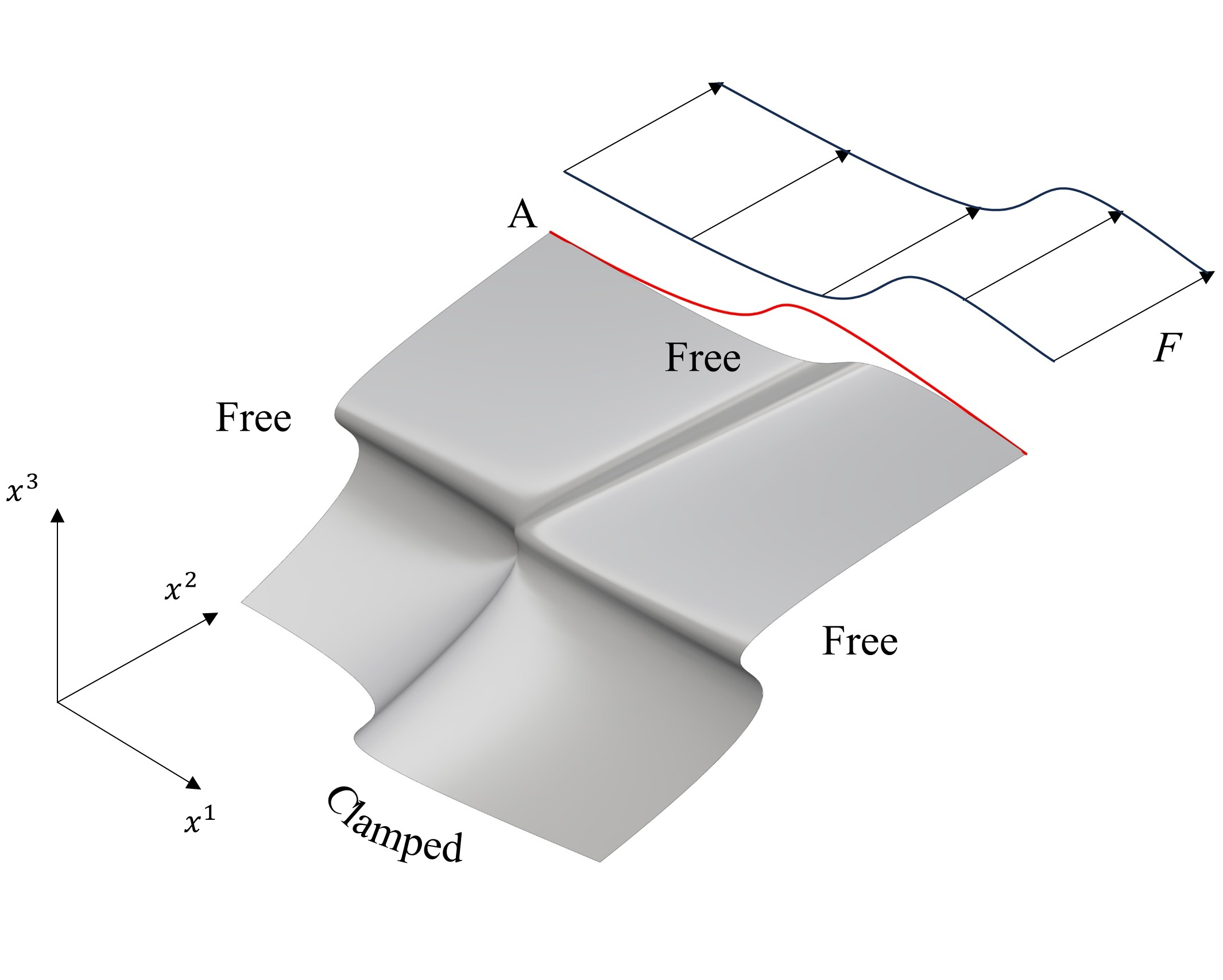} }
	\caption{\centering Schematic for the double curved free form NURBS based geometry: The red line shows the B-spline equivalent of the top side when the weights are set to 1}\label{fig_ffN_0}
\end{figure}

\begin{figure}[ht] 
	\center{\includegraphics[scale=0.9]{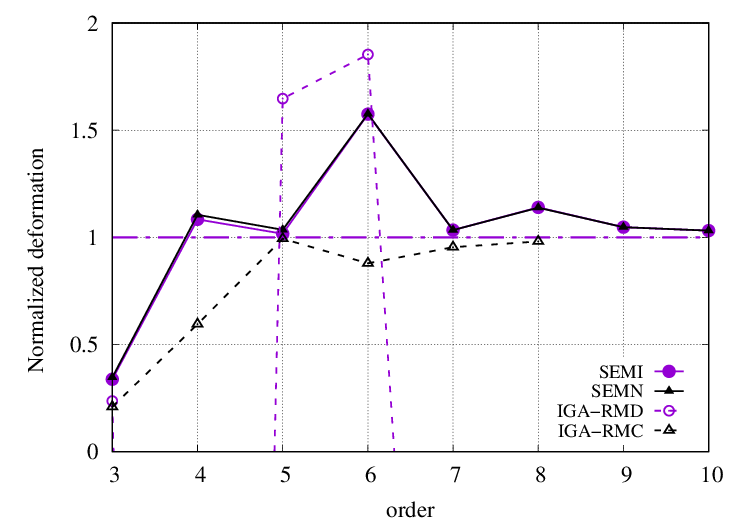} }
	\caption{\centering Double curved free form NURBS based geometry: Normalized deformation for \textit{p}-refinement  vs order of element}
	\label{fig_ffN_pref}
\end{figure}

\begin{figure}[ht]
\begin{subfigure}[t]{\linewidth}
	\includegraphics[width=120mm]{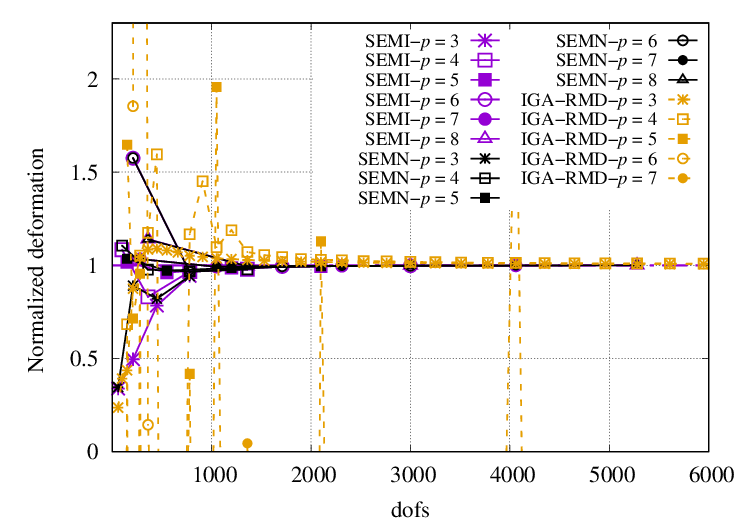} 
	\caption{comparison to IGA-RMD}
    \label{fig_ffN_href_a}
\end{subfigure}
\begin{subfigure}[t]{\linewidth}
	\includegraphics[width=120mm]{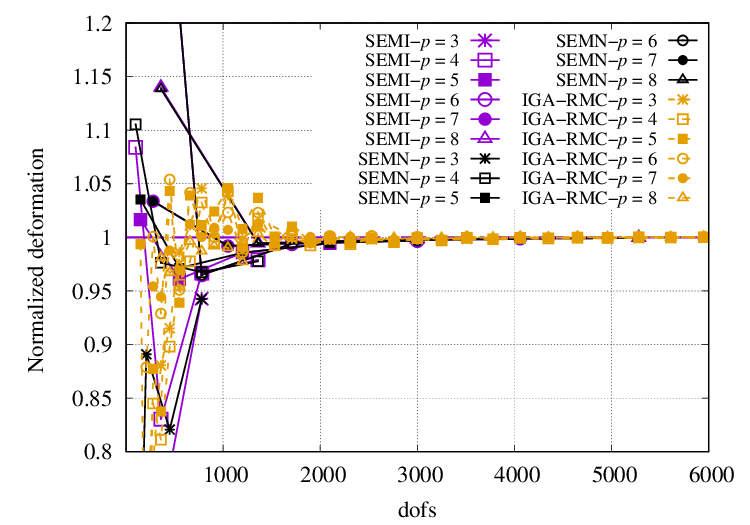}
	\caption{comparison to IGA-RMC}
    \label{fig_ffN_href_b}
\end{subfigure}
\caption{\centering Double curved free form NURBS based geometry: Normalized deformation vs DOFs for \textit{h}-refinement of SEMI and SEMN, compared to IGA}
\label{fig_ffN_href}
\end{figure}

\begin{figure}[ht]
        \begin{subfigure}[t]{\linewidth}
			\includegraphics[width=120mm]{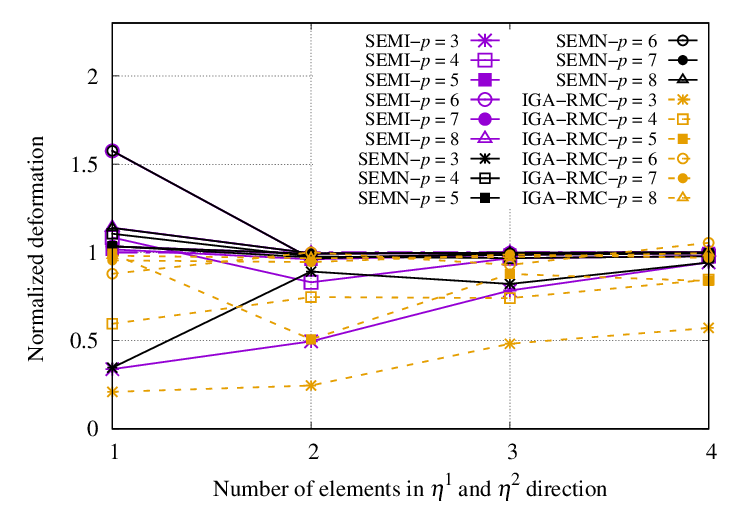}
			\caption{entire deformation range}\label{fig_ffN_href_elem_a}
	\end{subfigure}
    \begin{subfigure}[t]{\linewidth}
			\includegraphics[width=120mm]{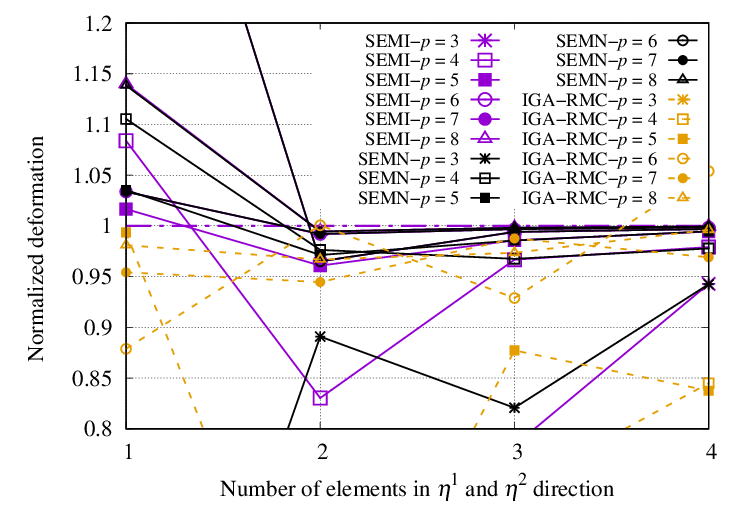}
			\caption{close up around the normalized deformation}\label{fig_ffN_href_elem_b}
	\end{subfigure}
	\caption{\centering Double curved free form NURBS based geometry: Normalized deformation vs number of elements  for \textit{h}-refinement of SEMI and SEMN, compared to IGA-RMC}
	\label{fig_ffN_href_elem}
\end{figure}

This example is the most sophisticated one as not only a double curved free form NURBS based shell surface, as shown in Fig.~\ref{fig_ffN_0}, is analyzed but also the geometry is adopted in a way that it reflects an abrupt change in the curvature. The thickness of the surface is 0.1 and the geometry of the surface can be constructed by using the knot vector and control points provided in \ref{app3}. A uniform line load of $F=10$ per unit length in $x^2$ direction is applied at the top edge of the surface. Young's modulus and Poisson's ratio are $1.2\cdot10^6$ and 0.3, respectively. The converged displacement at point A in direction of $x^2$ axis utilized to normalize the result is 0.3326847. Fig.~\ref{fig_ffN_pref} shows the results for $p$-refinement, using only one element. IGA-RMD does not result in accurate values at any order for this complicated structure. At lower orders, SEMI and SEMN outperform IGA-RMC. However, SEMI and SEMN exhibit oscillations even at higher orders, similar to IGA-RMC, though IGA-RMC has a smaller amplitude. Finally, SEMI and SEMN also converge to the correct value. The oscillations between the orders can be explained by the fact that while in uneven orders a node is present at the point of the strong change in curvature, this is not the case for even orders. Again, the exact representation of geometry by SEMN does not yield any gain of accuracy in comparison to SEMI. The comparison based on the number of DOFs is shown in Fig.~\ref{fig_ffN_href}. It confirms the difficulty of IGA-RMD in converging to the correct values at higher orders. At lower orders, only fine meshes provide good results in IGA-RMD, while at high orders there is divergence, see Fig.~\ref{fig_ffN_href_a}. On the other hand, unlike SEMI and SEMN, IGA-RMC exhibits osciallations with respect to DOFs regardless of the order of elements, though these fluctuations are alleviated in finer meshes, see Fig.~\ref{fig_ffN_href_b}. In Fig.~\ref{fig_ffN_href_elem}, the normalized deformation is plotted against the number of elements, using two different $y$-scales. From Fig. \ref{fig_ffN_href_elem_b}, it can be observed that the oscillations with respect to the number of elements are significantly more pronounced in IGA-RMC compared to SEMI and SEMN (for example, consider $p=5$). This figure clearly illustrates that while refining the mesh drastically increases the computational cost, it does not immediately mitigate the oscillations in IGA-RMC. In contrast to that, the oscillations in SEMI are significantly less pronounced.
\par
\begin{figure}[ht]
	\begin{subfigure}[t]{\linewidth}
			\includegraphics[width=120mm]{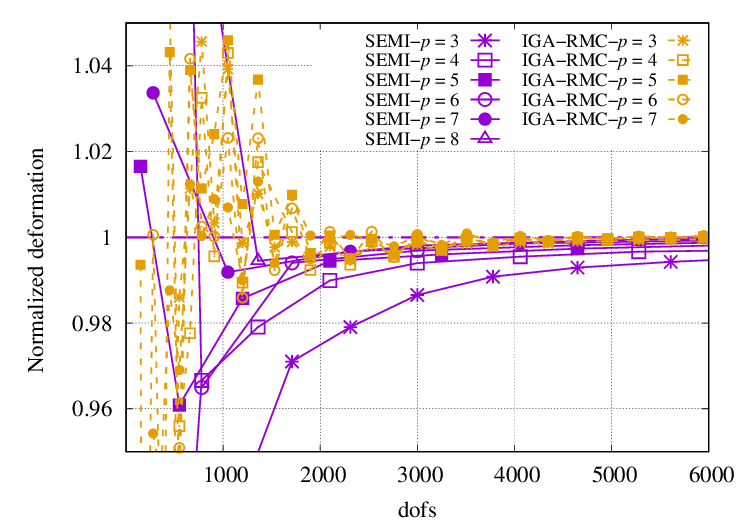} 
			\caption{$C^2$ continuity in IGA-RMC and $C^0$ continuity in SEMI }\label{fig_ffN_href_C0_a}
        \end{subfigure}
        \begin{subfigure}[t]{\linewidth}
			\includegraphics[width=120mm]{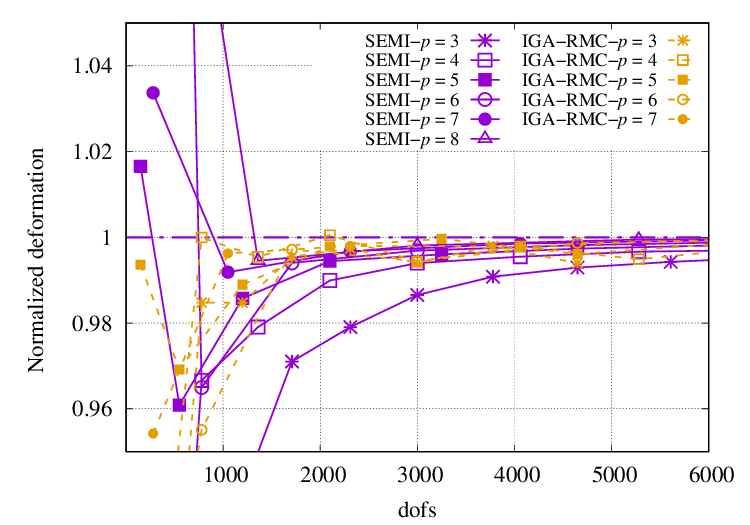}
			\caption{$C^0$ continuity, same as in SEMI}\label{fig_ffN_href_C0_b}
	\end{subfigure}
	\caption{\centering  Double curved free form NURBS based geometry: Normalized deformation vs DOFs for \textit{h}-refinement of SEMI, compared with IGA-RMC with different continuity of  elements at boundaries}
	\label{fig_ffN_href_C0}
\end{figure}

 In this paper, the term $h$-refinement inevitably refers to conventional $h$-refinement in SEM and $k$-refinement in IGA methods. However, this results in significantly fewer DOFs for the same number of elements in IGA compared to SEM. This discrepancy naturally raises curiosity about the impact of high continuity, which is the result of $k$-refinement, on IGA results. On the other hand, if only $C^0$ continuity is maintained in IGA, the influence of exact geometry definition can be directly compared with the results of SEMI. To explore this further, Fig.~\ref{fig_ffN_href_C0_a} and~\ref{fig_ffN_href_C0_b} are provided. In these figures, IGA-RMC results for order 8 are excluded due to instability in the solution. Fig.~\ref{fig_ffN_href_C0_a} is based on Fig.~\ref{fig_ffN_href}, but with a larger scale in the y-direction and the removal of SEMN results. To create the IGA-RMC results of Fig.~\ref{fig_ffN_href_C0_b}, knot insertion and order elevation is applied in a manner that at all element boundaries an inter-element continuity of $C^0$ prevails. Surprisingly, the downside of this continuity reduction in IGA-RMD is so significant that it can hardly converge to an accurate value at any order or with any number of elements. Therefore, it is excluded from the figures. The IGA-RMC results indicate that reducing continuity decreases the amplitude of oscillations in coarser meshes. However, the oscillations persist in finer meshes, and the results at large DOFs are not as accurate as those of IGA-RMC with higher continuity. It is to be noted, that the results of SEMI show proper convergence without oscillations for all orders, beginning from $2\times 2$ elements. For fine meshes and high orders, the accuracy of SEMI is higher than of IGA-RMC with $C^0$-continuity. This confirms that the SEM isoparametric representation of geometry is sufficient to achieve accurate and robust results. These figures also support the idea of using fewer higher-order SEMI elements rather than employing a fine mesh of low-order elements. It is to be noted, that especially for this most complicated example, SEMI performed as accurate as IGA-RMC, while being significantly more robust, i.e. showing less oscillations.

\section{Discussion and Conclusion} \label{sec_conclusion}
We have implemented a  Reissner--Mindlin shell formulation based on Rodrigues' rotation tensor for the additive update of the director vector within the spectral element method. In order to study the influence of the error entailed by the geometry approximation done in SEM methods, we compared the results with two different rotation interpolation  variants of an IGA-based Reissner-Mindlin shell formulation. The main findings can be summarized by:
\begin{itemize}
    \item The proposed SEMI formulation shows high accuracy and robustness. The employed discrete rotational formulation in combination with the Kronecker delta property of SEM prevents artificial thinning and has shown to be as accurate as the significantly more complicated continuous rotation formulation in IGA-RMC.
    \item In comparison to the more complicated SEMN formulation \cite{nima_2024}, where the geometry is exact, no significant disadvantage of SEMI could be found. Thus. the exact representation of geometry is not important in SEM.
    \item As discussed in detail in Sec.~\ref{sec_3}, in SEM, because the integration points coincide with the element nodes, in each integration point many components of the element stiffness matrix are equal to zero. Exploiting this characteristic feature significantly accelerates the generation of the element stiffness matrix in comparison to $p$-FEM/IGA formulations.
    \item The optimal strategy for SEM involves employing higher-order elements on coarser meshes, as opposed to utilizing a fine mesh composed of lower-order elements. In future work, static condensation of the DOFs related to the inner nodes of the elements \cite{payette_2014} will be tested to speed up computations.
    \item SEM methods can employ very high-order elements in large deformation analysis. The importance of this point is that according to \cite{payette_2014,nima_2024}, higher-order elements can significantly alleviate the effect of locking without requiring any additional treatment. A comparison in \cite{nima_2024} has shown that if the order is chosen high enough, the accuracy per DOF is comparable to the mixed method proposed in \cite{zou2020}, which relies on a much more complicated formulation.
     \item In the IGA formulations, instabilities arise at higher orders. IGA-RMC, due to its enhanced rotational formulation, exhibits greater stability than IGA-RMD. However, the maximum order of isogeometric elements, which provides stable results, is much less than the maximum stable order of spectral elements. These findings are in line with \cite{PITTON2018440}, which also proposes to use SEM functions to be able to make use high order basis functions. The deterioration of IGA results with rising order gets more pronounced, the more complicated the geometry is, i.e., high and highly changing curvature. The high continuity of IGA seems to be only beneficial for geometries with low changes in curvature.
    \item Regarding the small differences between the results of SEMI and IGA-RMC, it is clear that the exact geometry definition has a minimal effect on the accuracy of the results. In contrast to that, other choices such as the  rotational formulation have much more impact. Our results confirm the findings of~\cite{Oesterle.2022} that the real benefit of IGA is not the exact representation of geometry, but the high continuity. In proximity to the correct solution, the geometry approximation error is smaller than the interpolation approximation error. 
    \item Due to benefiting from  $k$-refinement, IGA yields higher accuracy per DOF, while SEMI yields higher accuracy per number of elements. If including efficient static condensation, this might change in favor of SEMI. The intended use of SEMI with static condensation while exploiting the cross pattern scheme of Sec.~\ref{sec_3} might yield a numerical scheme competitive with IGA. However, a fair evaluation and comparison of this requires the implementation within the same numerical framework, which we will tackle in upcoming works.
\end{itemize}
Basing on the aforementioned findings, in future works we aim at implementing a SEMI-based framework for the computation of highly complex CAD files, containing a multitude of trimmed NURBS patches. To keep meshing straightforward, all patches will be meshed independently, while the coupling will be performed using mortar methods~\cite{Bernardi.1993,Bernardi.1993,Dornisch.2017}, for which the Kronecker delta property is highly beneficial. This will further require the development of triangular Reissner--Mindlin SEM shell elements, such as proposed in~\cite{Petrovskiy_2016} for seven-parameter shells.

\section*{Declarations}

This work is part of a project funded by the Deutsche Forschungsgemeinschaft (DFG) under project number 503246947.

\appendix
\setcounter{table}{0}
\section{Construction of Local Nodal Coordinate System}
A local nodal Cartesian coordinate system at node $I$ can be constructed by the cross product of the director and the covariant tangent vector as below
\begin{align}
		&\mbn{A}_{3I}=\mbn{D} \\
		&\mbn{A}_{1I}=\frac{\mbn{G}^0_2 \times \mbn{D}}{\|\mbn{G}^0_2 \times \mbn{D}\|}\\
		&\mbn{A}_{2I}=\frac{\mbn{D} \times \mbn{A}_{1I}}{\|\mbn{D} \times \mbn{A}_{1I}\|}
\end{align}
\label{app1}
The nodal basis vectors $\mbn{A}_{1I}$ and $\mbn{A}_{2I}$ span a tangential plane to the shell midsurface in node $I$.

\section{Geometry of double curved free form surface}\label{app2}
\setcounter{table}{0}
 The double curved free form surface is a Coons patch uniquely defined by four boundary B-spline curves. The control points for these boundary curves are provided in Table~\ref{tab_ff_geom}. Since all control points have a weight value of 1, they collectively create a B-spline surface. The knot vector for all curves is
\begin{equation}
	\mbn{\Xi}=\left[0,0,0,0,\frac{1}{3},\frac{2}{3},1,1,1,1 \right]
\end{equation}
and the order of B-Spline basis functions is $p=3$ \cite{Dornisch.2013}.
\begin{table}[h!] 
   \caption{\centering Control points of four boundary curves of double curved free form surface }\label{tab_ff_geom}%
	\renewcommand{\arraystretch}{1.3}
	\begin{center}
			
			\begin{tabular}{@{}l|cccc@{}}
				\toprule
				 & Top & Bottom  & Left & Right\\
				\midrule
				Control  & 0,0,15    & 0,0,0   & 0,0,0  & 11,0,0  \\
				
				points &$\frac{11}{9}$,$\frac{2}{3}$,15 & 5,0,0  & 0,0,5 & 11,0,$\frac{8}{3}$ \\
				
				(x, y, z)&$\frac{11}{3}$,2,15 & 5,5,0  & 0,2,7 & 11,$\frac{2}{9}$,$\frac{62}{9}$ \\
				
				&$\frac{22}{3}$,4,15 & 10,5,0  & 0,2,10 & 11,$\frac{17}{9}$,$\frac{101}{9}$ \\
				
				&$\frac{88}{9}$,$\frac{16}{3}$,15 & 10,0,0  & 0,0,12 & 11,$\frac{13}{3}$,$\frac{41}{3}$ \\
				
				&11,6,15 & 11,0,0  & 0,0,15 & 11,6,15 \\
                \bottomrule
			\end{tabular}
	\end{center}
\end{table}

\section{Geometry of double curved free form NURBS based surface}\label{app3}
\setcounter{table}{0}
The double curved free form NURBS based surface is a Coons patch uniquely defined by three boundary B-splines and one boundary NURBS curves. The control points for these boundary curves are provided in table~\ref{tab_ffn_geom}.  The knot vector for all curves is
\begin{equation}
	\mbn{\Xi}=\left[0,0,0,0,1,1,1,1 \right]
\end{equation}
and the order of B-Spline and NURBS basis functions is $p=3$. Only two of the top points, among all the points, have weights different from 1, as shown in Table \ref{tab_ffn_geom}.
				
				
				


\begin{table}[h!] 
 \caption{\centering Control points of four boundary curves of double curved free form NURBS based surface}\label{tab_ffn_geom}%
	\renewcommand{\arraystretch}{1.3}
	\begin{center}
			
			\begin{tabular}{@{}l|ccccc@{}}
				\toprule
				 & Top & weights & Bottom  & Left & Right \\
				\midrule
				Control   & 0,4,5 & 1 & 0,0,0   & 0,0,0  & 5,4,0  \\
				
				points&4.972,7.188,5 & 1.5& 4.202,4.150,0  & 0,4.151,4.202 & 5,8.151,4.202 \\
				
				(x, y, z)&-1.303,5.255,5&0.5& 0.938, -1.546,0  & 0,-1.547,0.938 & 5,2.453,0.938 \\
				
				&5,8,5& 1& 5,4,0  & 0,4,5 & 5,8,5 \\
                \bottomrule
			\end{tabular}

	\end{center}
\end{table}



  \bibliographystyle{elsarticle-num} 
  \bibliography{references}

\begin{thebibliography}{10}
\expandafter\ifx\csname url\endcsname\relax
  \def\url#1{\texttt{#1}}\fi
\expandafter\ifx\csname urlprefix\endcsname\relax\def\urlprefix{URL }\fi
\expandafter\ifx\csname href\endcsname\relax
  \def\href#1#2{#2} \def\path#1{#1}\fi

\bibitem{hughes_2005}
T.~J.~R. Hughes, J.~A. Cottrell, Y.~Bazilevs, Isogeometric analysis: {{CAD}},
  finite elements, {{NURBS}}, exact geometry and mesh refinement, Computer
  Methods in Applied Mechanics and Engineering 194~(39) (2005) 4135--4195.
\newblock \href {https://doi.org/10.1016/j.cma.2004.10.008}
  {\path{doi:10.1016/j.cma.2004.10.008}}.

\bibitem{nguyen2015}
V.~Nguyen, C.~Anitescu, S.~Bordas, T.~Rabczuk, Isogeometric analysis: an
  overview and computer implementation aspects, Mathematics and Computers in
  Simulation 117 (2015) 89--116.
\newblock \href {https://doi.org/10.1016/j.matcom.2015.05.008}
  {\path{doi:10.1016/j.matcom.2015.05.008}}.

\bibitem{Piegl.1997}
L.~Piegl, W.~Tiller, {The NURBS book}, 2nd Edition, {Monographs in visual
  communications}, Springer, Berlin, 1997.

\bibitem{Kiendl.2009}
J.~Kiendl, K.-U. Bletzinger, J.~Linhard, R.~W{\"u}chner, {Isogeometric shell
  analysis with Kirchhoff-Love elements}, {Computer Methods in Applied
  Mechanics and Engineering} 198 (2009) 3902--3914.

\bibitem{Benson.2013}
D.~J. Benson, S.~Hartmann, Y.~Bazilevs, M.-C. Hsu, T.~J.~R. Hughes, {Blended
  isogeometric shells}, {Computer Methods in Applied Mechanics and Engineering}
  255 (2013) 133--146.
\newblock \href {https://doi.org/10.1016/j.cma.2012.11.020}
  {\path{doi:10.1016/j.cma.2012.11.020}}.

\bibitem{Benson.2010}
D.~J. Benson, Y.~Bazilevs, M.-C. Hsu, T.~J.~R. Hughes, {Isogeometric shell
  analysis: The Reissner--Mindlin shell}, {Computer Methods in Applied
  Mechanics and Engineering} 199 (2010) 276--289.
\newblock \href {https://doi.org/10.1016/j.cma.2009.05.011}
  {\path{doi:10.1016/j.cma.2009.05.011}}.

\bibitem{Oesterle.2017}
B.~Oesterle, R.~Sachse, E.~Ramm, M.~Bischoff, {Hierarchic isogeometric large
  rotation shell elements including linearized transverse shear
  parametrization}, {Computer Methods in Applied Mechanics and Engineering} 321
  (2017) 383--405.
\newblock \href {https://doi.org/10.1016/j.cma.2017.03.031}
  {\path{doi:10.1016/j.cma.2017.03.031}}.

\bibitem{zou2020}
Z.~Zou, M.~Scott, D.~Miao, M.~Bischoff, B.~Oesterle, W.~Dornisch, An
  isogeometric reissner–mindlin shell element based on bézier dual basis
  functions: Overcoming locking and improved coarse mesh accuracy, Computer
  Methods in Applied Mechanics and Engineering 370 (2020) 113283.
\newblock \href {https://doi.org/10.1016/j.cma.2020.113283}
  {\path{doi:10.1016/j.cma.2020.113283}}.

\bibitem{zou2022}
Z.~Zou, T.~Hughes, M.~Scott, D.~Miao, R.~Sauer, Efficient and robust
  quadratures for isogeometric analysis: Reduced gauss and gauss–greville
  rules, Computer Methods in Applied Mechanics and Engineering 392 (2022)
  114722.
\newblock \href {https://doi.org/10.1016/j.cma.2022.114722}
  {\path{doi:10.1016/j.cma.2022.114722}}.

\bibitem{ZOU2021}
Z.~Zou, T.~Hughes, M.~Scott, R.~Sauer, E.~Savitha, Galerkin formulations of
  isogeometric shell analysis: Alleviating locking with greville quadratures
  and higher-order elements, Computer Methods in Applied Mechanics and
  Engineering 380 (2021) 113757.
\newblock \href {https://doi.org/10.1016/j.cma.2021.113757}
  {\path{doi:10.1016/j.cma.2021.113757}}.

\bibitem{casquero2023}
H.~Casquero, K.~D. Mathews, Overcoming membrane locking in quadratic
  nurbs-based discretizations of linear kirchhoff–love shells: Cas elements,
  Computer Methods in Applied Mechanics and Engineering 417 (2023) 116523.
\newblock \href {https://doi.org/10.1016/j.cma.2023.116523}
  {\path{doi:10.1016/j.cma.2023.116523}}.

\bibitem{leonetti2024}
L.~Leonetti, D.~Magisano, G.~Garcea, Large rotation isogeometric shell model
  for alternating stiff/soft curved laminates including warping and interlayer
  thickness change, Computer Methods in Applied Mechanics and Engineering 424
  (2024) 116908.
\newblock \href {https://doi.org/10.1016/j.cma.2024.116908}
  {\path{doi:10.1016/j.cma.2024.116908}}.

\bibitem{Dornisch.2013}
W.~Dornisch, S.~Klinkel, B.~Simeon, {Isogeometric Reissner--Mindlin shell
  analysis with exactly calculated director vectors}, {Computer Methods in
  Applied Mechanics and Engineering} 253 (2013) 491--504.
\newblock \href {https://doi.org/10.1016/j.cma.2012.09.010}
  {\path{doi:10.1016/j.cma.2012.09.010}}.

\bibitem{Simo.1989b}
J.~C. Simo, D.~D. Fox, M.~S. Rifai, {On a stress resultant geometrically exact
  shell model. Part II: The linear theory; Computational aspects}, {Computer
  Methods in Applied Mechanics and Engineering} 73~(1) (1989) 53--92.

\bibitem{Dornisch.2015b}
W.~Dornisch, {Interpolation of Rotations and Coupling of Patches in
  Isogeometric Reissner--Mindlin Shell Analysis}, {Ph.D. thesis. Lehrstuhl
  f{\"u}r Baustatik und Baudynamik, RWTH Aachen}, 2015.

\bibitem{payette_2014}
G.~Payette, J.~Reddy, A seven-parameter spectral/hp finite element formulation
  for isotropic, laminated composite and functionally graded shell structures,
  Computer Methods in Applied Mechanics and Engineering 278 (2014) 664--704.
\newblock \href {https://doi.org/10.1016/j.cma.2014.06.021}
  {\path{doi:10.1016/j.cma.2014.06.021}}.

\bibitem{patera_1984}
A.~T. Patera, A spectral element method for fluid dynamics: {Laminar} flow in a
  channel expansion, Journal of Computational Physics 54~(3) (1984) 468--488.
\newblock \href {https://doi.org/10.1016/0021-9991(84)90128-1}
  {\path{doi:10.1016/0021-9991(84)90128-1}}.

\bibitem{poz_2005}
C.~Pozrikidis, Introduction {{To Finite And Spectral Element Methods Using
  Matlab}}, Chapman \& Hall/CRC, Boca Raton (FL), 2005.

\bibitem{PITTON2018440}
G.~Pitton, L.~Heltai, Nurbs-sem: A hybrid spectral element method on nurbs maps
  for the solution of elliptic pdes on surfaces, Computer Methods in Applied
  Mechanics and Engineering 338 (2018) 440--462.
\newblock \href {https://doi.org/10.1016/j.cma.2018.04.039}
  {\path{doi:10.1016/j.cma.2018.04.039}}.

\bibitem{nima_2024}
N.~Azizi, W.~Dornisch, A spectral finite element {{Reissner}}--{{Mindlin}}
  shell formulation with {{NURBS-based}} geometry definition, Computational
  Mechanics 74~(3) (2024) 537--559.
\newblock \href {https://doi.org/10.1007/s00466-024-02444-w}
  {\path{doi:10.1007/s00466-024-02444-w}}.

\bibitem{Gutierrez2016}
M.~{Gutierrez Rivera}, J.~Reddy, M.~Amabili, A new twelve-parameter spectral/hp
  shell finite element for large deformation analysis of composite shells,
  Composite Structures 151 (2016) 183--196.
\newblock \href {https://doi.org/10.1016/j.compstruct.2016.02.068}
  {\path{doi:10.1016/j.compstruct.2016.02.068}}.

\bibitem{zak2018}
A.~Żak, M.~Krawczuk, A higher order transversely deformable shell-type
  spectral finite element for dynamic analysis of isotropic structures, Finite
  Elements in Analysis and Design 142 (2018) 17--29.
\newblock \href {https://doi.org/10.1016/j.finel.2017.12.007}
  {\path{doi:10.1016/j.finel.2017.12.007}}.

\bibitem{Cohen.2007}
G.~Cohen, P.~Grob, Mixed higher order spectral finite elements for
  reissner–mindlin equations, SIAM Journal on Scientific Computing 29~(3)
  (2007) 986--1005.
\newblock \href {https://doi.org/10.1137/050642332}
  {\path{doi:10.1137/050642332}}.

\bibitem{Ambati.2025}
H.~Ambati, S.~Eisenträger, S.~Kapuria, Time-domain spectral bfs plate element
  with lobatto basis for wave propagation analysis, International Journal for
  Numerical Methods in Engineering 126~(3) (2025) e7617.
\newblock \href {https://doi.org/10.1002/nme.7617}
  {\path{doi:10.1002/nme.7617}}.

\bibitem{Simo.1989}
J.~C. Simo, D.~D. Fox, {On a stress resultant geometrically exact shell model.
  Part I: Formulation and optimal parametrization}, {Computer Methods in
  Applied Mechanics and Engineering} 72~(3) (1989) 267--304.

\bibitem{argyris_1982}
J.~Argyris, An excursion into large rotations, Computer Methods in Applied
  Mechanics and Engineering 32~(1) (1982) 85--155.
\newblock \href {https://doi.org/10.1016/0045-7825(82)90069-X}
  {\path{doi:10.1016/0045-7825(82)90069-X}}.

\bibitem{Dornisch.2016}
W.~Dornisch, R.~M{\"u}ller, S.~Klinkel, {An efficient and robust rotational
  formulation for isogeometric Reissner--Mindlin shell elements}, {Computer
  Methods in Applied Mechanics and Engineering} 303 (2016) 1--34.

\bibitem{wagner_2005}
W.~Wagner, F.~Gruttmann, A robust non-linear mixed hybrid quadrilateral shell
  element, International Journal for Numerical Methods in Engineering 64~(5)
  (2005) 635--666.
\newblock \href {https://doi.org/10.1002/nme.1387}
  {\path{doi:10.1002/nme.1387}}.

\bibitem{Bathe.2006}
K.-J. Bathe, {Finite Element Procedures}, 2nd Edition, Prentice Hall, Hoboken,
  2006.

\bibitem{grutt_2000}
F.~Gruttmann, R.~Sauer, W.~Wagner, Theory and numerics of three-dimensional
  beams with elastoplastic material behaviour, International Journal for
  Numerical Methods in Engineering 48~(12) (2000) 1675--1702.
\newblock \href
  {https://doi.org/10.1002/1097-0207(20000830)48:12<1675::AID-NME957>3.0.CO;2-6}
  {\path{doi:10.1002/1097-0207(20000830)48:12<1675::AID-NME957>3.0.CO;2-6}}.

\bibitem{Cottrell.2009}
J.~A. Cottrell, T.~J.~R. Hughes, Y.~Bazilevs, {Isogeometric analysis}, Wiley,
  Chichester, 2009.

\bibitem{dor_2016}
W.~Dornisch, R.~M{\"u}ller, S.~Klinkel, An efficient and robust rotational
  formulation for isogeometric {{Reissner}}--{{Mindlin}} shell elements,
  Computer Methods in Applied Mechanics and Engineering 303 (2016) 1--34.
\newblock \href {https://doi.org/10.1016/j.cma.2016.01.018}
  {\path{doi:10.1016/j.cma.2016.01.018}}.

\bibitem{belyt_1985}
T.~Belytschko, H.~Stolarski, W.~K. Liu, N.~Carpenter, J.~S. Ong, Stress
  projection for membrane and shear locking in shell finite elements, Computer
  Methods in Applied Mechanics and Engineering 51~(1) (1985) 221--258.
\newblock \href {https://doi.org/10.1016/0045-7825(85)90035-0}
  {\path{doi:10.1016/0045-7825(85)90035-0}}.

\bibitem{Oesterle.2022}
B.~Oesterle, F.~Geiger, D.~Forster, M.~Fröhlich, M.~Bischoff, A study on the
  approximation power of nurbs and the significance of exact geometry in
  isogeometric pre-buckling analyses of shells, Computer Methods in Applied
  Mechanics and Engineering 397 (2022) 115144.
\newblock \href {https://doi.org/10.1016/j.cma.2022.115144}
  {\path{doi:10.1016/j.cma.2022.115144}}.

\bibitem{Bernardi.1993}
C.~Bernardi, Y.~Maday, A.~T. Patera, {Domain Decomposition by the Mortar
  Element Method}, in: H.~G. Kaper, M.~Garbey, G.~W. Pieper (Eds.), {Asymptotic
  and Numerical Methods for Partial Differential Equations with Critical
  Parameters}, {Springer Netherlands}, Dordrecht, 1993, pp. 269--286.

\bibitem{Dornisch.2017}
W.~Dornisch, J.~St{\"o}ckler, R.~M{\"u}ller, {Dual and approximate dual basis
  functions for B-splines and NURBS -- Comparison and application for an
  efficient coupling of patches with the isogeometric mortar method}, {Computer
  Methods in Applied Mechanics and Engineering} 316 (2017) 449--496.
\newblock \href {https://doi.org/10.1016/j.cma.2016.07.038}
  {\path{doi:10.1016/j.cma.2016.07.038}}.

\bibitem{Petrovskiy_2016}
K.~A. Petrovskiy, A.~V. Vershinin, V.~A. Levin, Application of spectral
  elements method to calculation of stress-strain state of anisotropic
  laminated shells, IOP Conference Series: Materials Science and Engineering
  158~(1) (2016) 012077.
\newblock \href {https://doi.org/10.1088/1757-899X/158/1/012077}
  {\path{doi:10.1088/1757-899X/158/1/012077}}.

\end{thebibliography}





\end{document}